\documentclass[11pt,reqno]{amsart}
\usepackage{amsmath,amssymb,amsthm,mathtools}
\usepackage{enumitem}
\usepackage[hidelinks]{hyperref}
\newtheorem{theorem}{Theorem}[section]
\newtheorem{lemma}[theorem]{Lemma}
\newtheorem{proposition}[theorem]{Proposition}
\newtheorem{corollary}[theorem]{Corollary}
\theoremstyle{definition}
\newtheorem{example}[theorem]{Example}
\theoremstyle{remark}

\theoremstyle{plain}
\newtheorem*{maina}{Theorem A}
\newtheorem*{mainb}{Theorem B}
\newtheorem*{mainc}{Theorem C}
\newcommand{\K}{\mathrm{(K)}}
\newcommand{\m}{\mathfrak m}
\newcommand{\pd}{\operatorname{pd}}
\newcommand{\Hom}{\operatorname{Hom}}
\newcommand{\Ext}{\operatorname{Ext}}
\newcommand{\Spec}{\operatorname{Spec}}
\newcommand{\Sym}{\operatorname{Sym}}
\newcommand{\PP}{\mathbf P}
\title[Module-finite extensions over singular local rings]{Splitting module-finite extensions over singular local rings}
\author{J. D. V\'elez}
\address{Departamento de Matem\'aticas, Universidad Nacional de Colombia, Sede Medell\'in, Colombia}
\date{}
\hypersetup{pdftitle={Splitting module-finite extensions over singular local rings},pdfauthor={J. D. Velez}}
\begin{document}
\begin{abstract}
We study when a module-finite extension of a singular local ring splits as a module over the base. On a Gorenstein curve of positive residue characteristic, weak normality forces splitting for every extension of finite projective dimension. In residue characteristic two the converse holds when two is zero or a nonzerodivisor. A missing square root in the prescribed base gives an explicit nonsplit extension of generic rank four. We obtain a second curve criterion by gluing retractions on the branches and determine all evaluation ideals for coordinate crossings and for $x^2=y^2u z_1\cdots z_s$ in odd characteristic. For cubic extensions we prove a presentation theorem with arbitrarily many relations. It applies on the blow-up of an ideal generated by a regular sequence of length three, and on the blow-up of a three-dimensional local ring with a quadric tangent cone. A family of matrices satisfying the criterion exists over every residue field in the regular-sequence case.
\end{abstract}
\maketitle
\section{Introduction}
A \emph{module-finite extension} $R\subseteq S$ is an inclusion of commutative rings with the same identity such that $S$ is finitely generated as an $R$-module. It splits over $R$ if some $R$-linear functional $\rho:S\to R$ takes $1_S$ to $1$. All rings are Noetherian. The direct summand theorem settles the question over a regular base \cite{Andre}. For a singular local ring $(R,\m,k)$, Koh asked whether finite projective dimension of the extension suffices \cite{Koh}. We say that $R$ has the \emph{Koh property}, denoted by $\K$, if the following holds.
\[
 \text{Every module-finite extension }R\subseteq S
 \text{ with }\pd_R S<\infty\text{ splits over }R.
\]
We allow nilpotents and $R$-torsion in $S$.

The author's thesis and subsequent paper gave counterexamples to Koh's conjecture \cite{V93,Velez}. V\'elez and Fl\'orez later constructed lower-dimensional examples in characteristic two and in mixed characteristic \cite{VF}. Other positive results concern bigenerated extensions, bounded complexity, and annihilator and socle criteria \cite{GGV2016,VG2020}. Here we ask which curve singularities have $\K$, and which presentations force an individual extension to split in higher dimension.

\begin{maina}
Let $R$ be a one-dimensional Gorenstein local ring of residue characteristic $p>0$. If $R$ is weakly normal, it has $\K$. When $p=2$ and $2$ is zero or a nonzerodivisor in $R$, the converse holds. In equal characteristic $p$, weak normality is equivalent to $F$-purity.
\end{maina}

The point is that the values $\phi(1_S)$ of linear functionals form the conductor of an intermediate ring in the total quotient ring $Q(R)$. A determinant on a free resolution controls the powers of each scalar in that ring. Conversely, a missing square root yields a nonsplit extension over the given base, with minimal resolution $0\to R^6\to R^{10}\to S\to0$ (Theorem~\ref{realize}). The converse at odd primes remains undecided by this argument.

For an excellent one-dimensional local ring, the normalization of $(\widehat R)_{\mathrm{red}}$ is a product of complete discrete valuation rings. Say that its \emph{branch residue extensions are separable} if each of their residue fields is separable over $k$.

\begin{mainb}
An excellent seminormal local curve with separable branch residue extensions has $\K$. If $R$ is an excellent one-dimensional Cohen--Macaulay local ring with this separability property, residue characteristic two, and $2$ zero or a nonzerodivisor, then $\K$ is equivalent to seminormality. For an excellent Cohen--Macaulay local curve of equal characteristic $p>0$ with separable branch residue extensions, seminormality, weak normality and $F$-purity are equivalent.
\end{mainb}

Theorem B allows any number of branches and an imperfect residue field. Its proof starts with an exact criterion for a given extension of a fiber product of discrete valuation rings (Theorem~\ref{branch-split}). It uses the classical structure of seminormal and $F$-pure curves \cite{Swan,GT,GW,Gib,EH,HS}. Section~\ref{geometry} uses Cartier operators to pass from the generic point of a crossing to the whole singularity.

The last result in this introduction concerns cubic extensions. Let $I\subseteq\m$ and let $U$ be a $(b+2)\times b$ matrix with entries in $I$. On the closed fiber of the blow-up of $I$, its entries induce linear forms from $I/\m I$. We denote the resulting matrix by $\overline U_I$. Full column rank means rank $b$ at every geometric point of this fiber.

\begin{mainc}
Let $R\subseteq S$ be a module-finite extension with $S/R\simeq\operatorname{coker}U$, where $b\ge1$. Suppose $\overline U_I$ has full column rank and one of the following holds.
\begin{enumerate}[label=\textup{(\roman*)}]
\item $I$ is generated by an $R$-regular sequence of length three.
\item $I=\m$ and
$\operatorname{gr}_{\m}R\simeq k[X_0,X_1,X_2,X_3]/(f_2)$ for a nonzero quadratic form $f_2$.
\end{enumerate}
Then $R\subseteq S$ splits. In either case arbitrary changes of the entries in $\m I$ preserve the criterion.
\end{mainc}

There is no restriction on $b$, on the residue characteristic, or on the rank of $f_2$. For $I=(a_0,a_1,a_2)$ as in (i), multiplication of polynomials by $a_0+a_1Z+a_2Z^2$ gives an explicit admissible matrix for every $b$ (Corollary~\ref{toeplitz}). The proof uses the binary cubic of a rank-three algebra \cite[Theorem~2.1]{Wood}. In residue characteristic three its restriction to the exceptional fiber is controlled by a two-step filtration of a symmetric cube. For a hypersurface tangent cone of higher degree, Section~\ref{higher-surfaces} gives separate criteria for two and three relations.

In particular, Corollary~\ref{three-entries} handles all degrees and all numbers of relations when the entry ideal has three generators and the linear matrix has full rank on the exceptional surface.

Sections~\ref{homol} and \ref{ranks} treat individual extensions by means of $\Ext^1$, the tangent cone, and multiplicity. The quadratic result refines the argument of Hochster--McLaughlin \cite[Theorem~2.1]{HM}, replacing local factoriality by triviality of the reflexive determinant of the specified algebra. The target may have zero divisors.

\subsection*{Notation}
Write $M^*=\Hom_R(M,R)$, $\ell_R(M)$ for length, and
$H^0_{\m}(M)=\bigcup_n(0:_M\m^n)$. The \emph{evaluation ideal} of an extension is
\[
 J(S/R)=\{\phi(1_S):\phi\in\Hom_R(S,R)\}.
\]
It is the trace ideal, since $\phi(s)$ is the value at $1$ of $t\mapsto\phi(st)$. Thus $J(S/R)=R$ exactly when the extension splits. Evaluation commutes with flat base change by finite presentation, so a faithfully flat extension of the base detects splitting. Also, $1_S\notin\m S$ by Nakayama's lemma. If $S$ is free, its unit extends to a basis. In particular every local ring of depth zero has $\K$ by Auslander--Buchsbaum \cite{BH}.

A seminormal (respectively weakly normal) ring is reduced and admits no proper module-finite extension inside $Q(R)$ inducing a bijection on spectra and isomorphisms (respectively purely inseparable extensions) on residue fields. We use the square--cube criterion for seminormality and the root criterion of Lemma~\ref{wn}. In characteristic $p$, $F_*^eM$ denotes restriction of scalars along the $e$th Frobenius. The ring is $F$-pure when $R\to F_*R$ remains injective after tensoring with any module. We impose no finiteness condition on Frobenius.

One rank observation applies in every dimension. For a module-finite extension of finite projective dimension, lift multiplication by $s\in S$ to a chain map $f_\bullet(s)$ on a finite free resolution $P_\bullet$ of $S$. The alternating trace
\begin{equation}\label{perfect-trace}
 \tau(s)=\sum_i(-1)^i\operatorname{tr}(f_i(s))
 \quad\text{satisfies}\quad
 \tau(1)=r=\sum_i(-1)^i\operatorname{rank}P_i .
\end{equation}
Two lifts are homotopic. The alternating trace of $dh+hd$ vanishes by $\operatorname{tr}(AB)=\operatorname{tr}(BA)$, so $\tau$ is well defined and $R$-linear. At each minimal prime $S$ is free by Auslander--Buchsbaum, of rank $r$. Thus $r\in J(S/R)$ and $\tau/r$ is a retraction whenever $r$ is a unit. In particular all finite-projective-dimension extensions of generic rank three split if $3$ is a unit.

\section{Gorenstein curves and the evaluation ideal}\label{scalar}
Let $R$ be a one-dimensional Gorenstein local ring. Fix a regular parameter $x$. Every nonunit nonzerodivisor generates an $\m$-primary ideal and hence divides a power of $x$. Thus $Q=Q(R)=R[1/x]$. For any module-finite extension $R\subseteq S$, set
\[
 T=H^0_{\m}(S),\qquad U=S/T,\qquad B=U\cap Q\ \subseteq S[1/x].
\]
Here $T$ consists of the elements killed by a power of $x$, so it is an ideal and meets $R$ trivially. The ring $U$ embeds in $S[1/x]$, and $B$ is a finite birational $R$-algebra. Its conductor $(R:B)$ is the ideal of all $a\in R$ such that $aB\subseteq R$.
In the intersection defining $B$, $Q$ denotes the image of the
scalar map $Q\hookrightarrow S[1/x]$.  Finiteness follows
because $B$ is an $R$-submodule of the finite module $U$.
\begin{proposition}\label{eval}
The ideal $\{\phi(1):\phi\in\Hom_R(S,R)\}$ equals $(R:B)$. In particular, the extension splits precisely when $B=R$, and $B/R=H^0_{\m}(U/R)$.
\end{proposition}
\begin{proof}
Maps to $R$ kill $T$ by regularity of $x$. The class of $u\in U$ in $U/R$ is killed by a power of $x$ exactly when $u$ belongs to $Q$. Likewise, if $x u\in B$ then $u\in B$, so $U/B$, when nonzero, has depth one.

The point where the Gorenstein hypothesis enters is extension of
functionals across $B\subseteq U$.

For any finite module $M$ on which $x$ is regular, $\Ext^1_R(M,R)=0$. Indeed $G=R/xR$ is Artinian Gorenstein, hence self-injective \cite{BH}. A free resolution of $M$ stays exact modulo $x$, since $\operatorname{Tor}_{>0}^R(M,G)=0$. Applying $\Hom_G(-,G)$ gives $\Ext^1_R(M,G)=0$. The coefficient sequence for multiplication by $x$ makes $x$ surjective on $\Ext^1_R(M,R)$, which vanishes by Nakayama. Applying this to $U/B$ shows that every functional on $B$ extends to $U$. After inverting $x$, a functional on $B$ becomes a $Q$-linear map $Q\to Q$, hence multiplication by its value at $1$. Its possible values are precisely $(R:B)$.
\end{proof}

\subsection{Generic rank and weak normality}\label{roots}
The scalar obstruction is constrained by generic rank.  A characteristic
polynomial computed on a finite free resolution supplies the necessary
integrality relations in both equal and mixed characteristic.
\begin{proposition}\label{rank}
Assume $\pd_R S<\infty$ and $\operatorname{char} k=p>0$. Then $S_Q$ is free over $Q$ of rank $r\ge1$. For $b\in B$ one has
\[
 \binom rj b^j\in R\quad(0\le j\le r),\qquad
 p^e b\in R,\qquad b^{p^e}\in R,\qquad e=v_p(r).
\]
\end{proposition}
\begin{proof}
By Auslander--Buchsbaum, $S$ admits a free resolution $0\to F_1\to F_0\to S\to0$. Over each local factor of the Artinian ring $Q$, the localized module is free by the same formula and is nonzero by scalar injectivity. Its rank is the common integer $r=\operatorname{rank} F_0-\operatorname{rank} F_1$.

Lift $b$ to $s\in S$ and lift multiplication by $s$ to endomorphisms $f_i$ of the resolution. The series $\det(1-Zf_0)/\det(1-Zf_1)$ belongs to $R[[Z]]$. Over $Q$ the resolution splits. In an adapted basis $f_0$ is block triangular with diagonal blocks $f_1$ and $bI_r$. The series therefore equals $(1-Zb)^r$ in $Q[[Z]]$. Comparing coefficients gives the first assertion. Write $r=p^e m$ with $p\nmid m$. Since $m$ is a unit of $R$, the linear coefficient gives $p^e b\in R$. In $\mathbf F_p[Z]$ the coefficient of $Z^{p^e}$ in $(1+Z)^r=(1+Z^{p^e})^m$ is $m$, so $\binom r{p^e}$ is also a unit in $R$.
\end{proof}
For $p\nmid r$ this gives $B=R$, recovering the splitting in \eqref{perfect-trace}. 

\begin{lemma}\label{wn}
Let $R$ be one-dimensional Cohen--Macaulay local of residue characteristic $p>0$. The following conditions are equivalent.
\begin{enumerate}
\item $R$ is weakly normal.
\item $q\in Q(R)$ and $q^p\in R$ imply $q\in R$.
\item $q\in Q(R)$ and $pq,q^p\in R$ imply $q\in R$.
\end{enumerate}
\end{lemma}
\begin{proof}
For $q^p\in R$ the extension $R\subseteq R[q]$ is module-finite and is an isomorphism at minimal primes, where a regular denominator for $q$ is invertible. The closed fiber is a nonzero quotient of $k[t]/(t^p-\overline{q^p})$, which has a unique prime and purely inseparable residue field. A one-dimensional local spectrum has only its minimal primes and its closed point. Thus the extension is weakly subintegral, proving (1)$\Rightarrow$(2).

Clearly (2) implies (3).  To prove the reverse implication, start with
$q^p\in R$.  The finite module $R[q]/R$ vanishes at the minimal
primes, so it has finite length.  Since $p\in\m$, there is an
$N\ge0$ with $p^Nq\in R$.  If $N>0$, set $q'=p^{N-1}q$.
Then $pq'\in R$ and $(q')^p=p^{p(N-1)}q^p\in R$, so (3)
gives $q'\in R$.  Descending induction gives $q\in R$.

Assume (2) and choose a regular parameter $x$ of $R$. A nilpotent $z\in R$ has $z/x^n\in R$ for every $n$, by repeated root closure applied to a vanishing $p$-power. Thus $z\in\bigcap_n x^nR=0$ by Krull intersection, and $R$ is reduced.
Let $C\subset Q$ be finite weakly subintegral over $R$. It is local, say with maximal ideal $\mathfrak n$, and $C/R$ has finite length. The conductor $I=(R:C)$ contains a power of $\m$ and is an ideal of $C$, so contains a power of $\mathfrak n$ since $\sqrt{\m C}=\mathfrak n$. We claim that $I$ is radical in $C$. If $b^N\in I$, then for any $c\in C$ and $p^h\ge N$ we have $(bc)^{p^h}\in I\subseteq R$. Repeated use of (2) gives $bc\in R$ for every $c$, proving $b\in I$. If $C\ne R$, then $I$ is proper and therefore $I=\mathfrak n$. For $b\in C$ some $p$-power of its residue belongs to $k$ by pure inseparability. Lifting that residue to $R$ gives $b^{p^h}\in R+\mathfrak n=R$, so $b\in R$ after all.
\end{proof}
\begin{proof}[Proof of the sufficient direction of Theorem A]
By Proposition~\ref{rank}, each element of $B$ has a $p$-power in $R$. Lemma~\ref{wn} gives $B=R$ by iteration, and Proposition~\ref{eval} gives the retraction.
\end{proof}

\begin{corollary}\label{fpure}
If $R$ is one-dimensional Gorenstein of characteristic $p>0$, then weak normality is equivalent to $F$-purity. In particular in characteristic two, $\K$ is also equivalent to $F$-purity. No finiteness condition on Frobenius is required.
\end{corollary}
\begin{proof}
Write $F_*M$ for restriction of scalars along Frobenius. By local duality for a Gorenstein ring \cite{BH}, the injective hull $E$ of $k$ identifies with $H^1_{\m}(R)=Q/R$. Localizing the twisted module gives $Q\otimes_R F_*R\simeq F_*Q$, via $v\otimes F_*a\mapsto F_*(v^p a)$. Taking cokernels of the maps from $F_*R$ therefore identifies $E\otimes_R F_*R$ with $F_*(Q/R)$. Under this identification the map obtained by tensoring Frobenius with $E$ sends $v+R$ to $F_*(v^p+R)$. Its injectivity is precisely root closure.

For completeness, a map $i:R\to M$ is pure if $E\otimes i$ is injective. Indeed the latter injection splits by injectivity of $E$. Adjunction sends a left inverse to an $R$-linear map $M\to\operatorname{End}_R(E)=\widehat R$ whose composite with $i$ is the completion map. After tensoring with any module the latter is still injective by faithful flatness. This proves purity. The converse follows from its definition. The result follows from Lemma~\ref{wn} and Theorem A.
\end{proof}

\section{Realizing a missing square root}\label{construction}
For the converse to Theorem A we construct an algebra over the given ring from a fraction in its total quotient ring.
\begin{theorem}[Realization over the prescribed base]\label{realize}
Let $(R,\m,k)$ have residue characteristic two. Suppose that either
$2=0$ and $q\in Q(R)\setminus R$ satisfies $q^2\in R$, or
$2$ is a nonzerodivisor and $q\in Q(R)\setminus R$ satisfies $2q,q^2\in R$.
There is a nonsplit module-finite extension $R\subseteq S$ with $S$ local,
$S\otimes_R Q(R)\simeq Q(R)^4$ as a module, and minimal resolution
$0\to R^6\to R^{10}\to S\to0$.
\end{theorem}
\begin{lemma}\label{model}
Take either $P=\mathbf F_2[X,(1+X)^{-1},Y]/(Y^2)$ with $x=X$, or $P=\mathbf Z_{(2)}[Y]/(Y^2)$ with $x=2$. There is a commutative $P$-algebra generated as a module by $1,u,v,w,d$, with relation matrix
\[
 D=\begin{pmatrix}0&x^3Y&0\\0&x^4&0\\Y&0&0\\x^3&Y&0\\0&0&x^2\end{pmatrix}.
\]
It satisfies $x^6(Y+xu)=0$. For any specialization to a local ring $A$ with $x\mapsto a$ regular in its maximal ideal and $Y\mapsto\eta$, the specialized matrix is injective, $A$ embeds in the algebra, and the closed fiber is
\[
 k_A[u,v,w]/(u^3,v^2,w^2,uv,uw,vw).
\]
\end{lemma}
\begin{proof}
Write $M=\operatorname{coker}D$ and $e_0$ for its first generator. The operators $U,V,W$ of Appendix~\ref{cert} induce commuting endomorphisms of $M$. Their values at $e_0$ are $u,v,w$. In the notation of the appendix,
\[
 (U^2-(x^3+h)U-x^2Y\,I)e_0=(1+x)d.
\]
As $1+x$ is invertible in both coefficient rings, $e_0$ is cyclic for the commutative algebra $C=P[U,V,W]$. An operator in $C$ killing $e_0$ kills its entire $C$-orbit, so evaluation at $e_0$ is injective. Hence $C\simeq M$ by evaluation, giving $M$ a unital commutative algebra structure.

The first relation implies $x^3Yw=0$ after multiplication by $Y$, and the second then gives $x^6(Y+xu)=0$. Tensoring with $A$ preserves the presentation. If a linear combination of the three columns has all its entries zero except possibly the first, rows two, four and five force all three coefficients to vanish, in that order, by regularity of $a$. This gives exactness and scalar injectivity. In the closed fiber the matrix vanishes. The operators there give $u^2=d$ and every other product of two generators in $u,v,w,d$ is zero (for products with $d$ use $d=u^2$). This yields the stated local closed fiber.
\end{proof}
\begin{proof}[Proof of Theorem~\ref{realize}]
In equal characteristic write $q=b/a$ with $a$ regular and $b\notin aR$. Set $c=q^2$. Both $a,b$ lie in $\m$, the latter because $b^2=a^2c$. Form
\[
 A=R[t]/(t^2-c),\qquad \eta=b+at .
\]
This is free on $1,t$, and $\eta^2=0$. In mixed characteristic put $a=2$, $b=2q$, $c=q^2$. Now $b^2=4c$, $b\in\m$ and $b\notin 2R$ by regularity of $2$. Put
\[
 A=R[t]/(t^2-bt+c),\qquad \eta=2t-b .
\]
Again $A$ is free on $1,t$ and $\eta^2=b^2-4c=0$. In either case the closed fiber of $A$ is $k[t]/(t^2-\bar c)$, which has one prime. Thus $A$ is local. The element $a$ stays regular and $\eta$ lies in its maximal ideal.

Specialize Lemma~\ref{model} over $A$ and call the resulting algebra $S$. It is finite and local, since its closed fiber over $A$ is local. The maps $R\hookrightarrow A\hookrightarrow S$ are injective. Restriction of the resolution $0\to A^3\to A^5\to S\to0$ gives the indicated ranks over $R$. To see minimality, the only block needing inspection is multiplication by $\eta$ on $A$, whose matrices in the two cases are respectively
\[
 \begin{pmatrix}b&ac\\a&b\end{pmatrix},\qquad
 \begin{pmatrix}-b&-2c\\2&b\end{pmatrix}.
\]
Thus all entries belong to $\m$. Inverting $a$ eliminates $w,u,d$ successively, leaving a free module on $1,v$ over $A[1/a]$. This proves the generic rank assertion.

Finally $a^6(\eta+au)=0$. Every $R$-linear map to $R$ kills $\eta+au$. For a retraction this would give $b+a\rho(t+u)=0$ in equal characteristic, and $2\rho(t+u)-b=0$ in mixed characteristic. Both imply $b\in aR$, a contradiction. Notice that the obstruction is tested over $R$ even though the presentation was constructed over $A$.
\end{proof}
\begin{proof}[Necessity in Theorem A]
If weak normality fails, Lemma~\ref{wn}(3) fails for $p=2$. The fraction witnessing failure satisfies the hypotheses of Theorem~\ref{realize} in each of the two cases of the converse.
\end{proof}
\begin{example}[A cusp and a mixed-characteristic domain]\label{root-examples}
Let $k$ have characteristic two and let
$R=k[[t^2,t^3]]$.  This is the one-dimensional Gorenstein domain
$k[[x,y]]/(y^2-x^3)$.  The fraction $q=t=y/x$ does not belong to
$R$, while $q^2=t^2\in R$.  Theorem~\ref{realize} therefore constructs a
nonsplit rank-four algebra of projective dimension one over this
specific cusp.

For the mixed-characteristic example put
$R=\mathbf Z_2[z]/(z^2-8)$.  It is finite free on $1,z$ over
$\mathbf Z_2$, and its reduction modulo two is
$\mathbf F_2[z]/(z^2)$, so it is complete and local.  The
polynomial is irreducible over $\mathbf Q_2$, since a square in
$\mathbf Q_2$ has even valuation whereas $v_2(8)=3$.  Thus $R$
is a domain.  The element $q=z/2$ lies in its fraction field,
has $q^2=2$ and $2q=z$, and is not in $R$ by the free basis
$1,z$.  The mixed-characteristic part of Theorem~\ref{realize} applies.
\end{example}
\begin{corollary}
A nonreduced Noetherian local ring of characteristic two satisfies $\K$ exactly when its depth is zero. If two is regular and the residue characteristic is two, $\K$ likewise forces reducedness, in every dimension. In each positive-depth case a nonreduced base admits a nonsplit module-finite local extension of generic rank two and Betti numbers $(5,3)$.
\end{corollary}
\begin{proof}
At depth zero, finite projective dimension forces freeness by Auslander--Buchsbaum. In a nonzero finite free algebra the unit is outside $\m S$ by Nakayama, and therefore extends to a basis. At positive depth choose $0\ne z\in R$ with $z^2=0$. In characteristic two take $a$ to be a sufficiently large power of a regular element of $\m$ so that $z\notin aR$, using Krull intersection, and put $\eta=z$. When two is regular take $a=2$ and write $z=2^n\eta$ for $n$ maximal. Then $\eta^2=0$ and $\eta\notin 2R$. Apply Lemma~\ref{model} directly over $R$. Its relation $a^6(\eta+au)=0$ precludes a retraction, and the same elimination and minimality arguments give the assertions.
\end{proof}

\section{Retractions from the branches}\label{branches}
Let $V_i$ ($1\le i\le s$) be DVRs, with maximal ideals $\mathfrak n_i=(\pi_i)$ and identified residue field $k$. Let $R=\prod_k V_i$ be the subring of tuples with common residue in $k$. Its maximal ideal is $J=\prod_i\mathfrak n_i$. The coordinate maps are onto, each summand of $J$ is Noetherian over its branch, and consequently $R$ is Noetherian. Its normalization is $\prod_i V_i$, since the quotient by $R$ is a finite-dimensional $k$-space and inverting $(\pi_i)_i$ gives the product of the fraction fields. Thus $R$ is reduced of dimension one, and its minimal primes are the coordinate kernels $P_i$.
\begin{theorem}[Branch criterion]\label{branch-split}
The ring $R=\prod_k V_i$ has $\K$ in every characteristic. When the branches are complete, a module-finite extension $R\subseteq S$ splits exactly when some local factor of $S$ is nonzero at every minimal prime of $R$.
\end{theorem}
\begin{lemma}\label{glue}
Every module-finite extension $R\subseteq T$ with $T$ local admits a retraction.
\end{lemma}
\begin{proof}
Put $T_i=T/P_iT$. Then $V_i\to T_i$ is injective. Indeed the tuple supported on branch $i$ with value $\pi_i$ annihilates $P_iT$. If $r\in R\cap P_iT$, multiplication by that tuple gives $\pi_i r_i=0$ in $V_i$. Thus $r\in P_i$. The $V_i$-torsion in $T_i$ is a proper ideal $N_i$, hence is contained in the unique maximal ideal of $T_i$. The closed fibers $T_i/\mathfrak n_iT_i$ all identify with $E=T/JT$.

Let $L$ be the residue field of $T$ and choose a $k$-linear functional $L\to k$ taking $1$ to $1$. The resulting functional $\lambda:E\to k$ kills the image of each $N_i$. Since $T_i/N_i$ is finite free over $V_i$, the induced functional on its closed fiber lifts, by choosing a basis, to a $V_i$-linear functional. Its value at $1$ is congruent to $1$ modulo $\mathfrak n_i$. Divide by this unit and compose with $T\to T_i/N_i$ to get $\psi_i:T\to V_i$. All the $\psi_i$ reduce to $\lambda$ and take $1$ to $1$. Thus $(\psi_i)_i$ takes values in $R$ and is the desired $R$-linear retraction.
\end{proof}
\begin{proof}[Proof of Theorem~\ref{branch-split}]
If the $V_i$ are complete then $R$ is complete, since $J^n=\prod_i\mathfrak n_i^n$. A finite $R$-algebra $S$ is therefore a product of nonzero local algebras $T_j$. A factor is nonzero on every generic branch exactly when the map $R\to T_j$ has zero kernel. Indeed, the generic localizations are fields and $R$ embeds in their product. Such a factor gives a retraction by Lemma~\ref{glue}. Conversely, decomposing a retraction as a sum of functionals on the factors gives $1=\sum_j\rho_j(1_{T_j})$. Some summand is a unit in the local ring $R$, forcing the corresponding factor map to be injective.

If $S$ has finite projective dimension, so does every $T_j$. The alternating sum of the ranks in a finite free resolution of $T_j$ equals its vector-space dimension on each generic branch. Some factor has positive rank on the first branch since $R\subseteq S$. It therefore has positive rank on every branch and yields the required retraction.

In general $\widehat R=\prod_k\widehat V_i$ by the same computation of $J^n$. Flat base change preserves the extension and its finite free resolution. For any flat extension $R\to R'$ the evaluation ideal extends to the evaluation ideal over $R'$, since $S$ is finitely presented and $\Hom_R(S,R)$ commutes with flat base change. Faithful flatness of completion thus descends the retraction.
\end{proof}

\subsection{Separable residue fields and intrinsic curves}

Let $(R,\m,k)$ be an excellent one-dimensional local ring and put
$A=(\widehat R)_{\mathrm{red}}$.  Its normalization is finite and has
the form $\overline A=\prod_iV_i$, where the $V_i$ are complete DVRs.
Their residue fields $k_i$ are finite over $k$.  We use the condition
\[
 \tag{SB}\label{separable-branches}
 \text{Every residue field extension }k_i/k\text{ is separable.}
\]
The use of the reduced completion makes this condition meaningful
before reducedness or seminormality of $R$ is known.  It is automatic
when $k$ is perfect, but the argument only needs
\eqref{separable-branches}.

We first recall the classical conductor description of a complete
seminormal curve \cite{Swan,EH}.

\begin{lemma}\label{conductor-description}
Let $R$ be a complete one-dimensional seminormal local ring with finite
normalization $V=\prod_iV_i$, where the $V_i$ are complete DVRs.  Put
$J=\prod_i\mathfrak n_i$ and $k_i=V_i/\mathfrak n_i$.  Then $J=\m$
as subsets of $V$, and $R$ is the inverse image of the diagonal image
of $k$ under the residue map $V\to\prod_i k_i$.
\end{lemma}
\begin{proof}
The finite module $V/R$ is supported at the closed point.  Its
annihilator, the conductor $I=(R:V)$, therefore contains a power of
$\m$.  Since $I$ is also an ideal of $V$ and $\sqrt{\m V}=J$, it
contains a power of $J$.

For $b\in J$, all sufficiently high powers of $b$ lie in $I\subseteq R$.
The square--cube criterion for seminormality \cite{Swan} says that
$c^2,c^3\in R$, for $c\in Q(R)$, imply $c\in R$.  If all powers $b^n$
with $n\geq N$ belong to $R$, applying this criterion to $b^n$ gives
the same conclusion for $n\geq\lceil N/2\rceil$.  Iteration yields
$b\in R$, so $J\subseteq R$.  The quotient $R/J$ is a reduced Artinian
local subring of $\prod_i k_i$, hence a field.  It is therefore the
residue field $k$ of $R$, proving both assertions.
\end{proof}

\begin{proof}[Proof of Theorem B]
Suppose first that $R$ is seminormal and satisfies
\eqref{separable-branches}.  Completion preserves seminormality for
excellent local rings \cite[Cor.~5.3]{GT}.  The evaluation-ideal
argument in the proof of Theorem~\ref{branch-split} permits descent, so we may assume
that $R$ is complete.  Let $V=\prod_iV_i$ be its normalization.
Lemma~\ref{conductor-description} identifies $R$ with the inverse image
of $k\subseteq\prod_i k_i$.

By \eqref{separable-branches}, a finite Galois extension $l/k$ splits
all the fields $k_i$.  Lift the monic minimal polynomial of a primitive
element of $l/k$ to a monic polynomial $f\in R[t]$.  Then
$R'=R[t]/(f)$ is finite free, local and complete, with residue field
$l$.  Its derivative is a unit at the class of $t$, as can be checked
in $l$, so $R'/R$ is finite \'etale.

Put $V'=V\otimes_RR'$.  On each branch, the field factors of
\[
 k_i\otimes_k l\simeq l^{[k_i:k]}
\]
lift by henselianity to local factors of the finite \'etale
$V_i$-algebra $V_i\otimes_RR'$.  Each factor is a complete DVR because it is
flat over $V_i$, its quotient by $\pi_i$ is $l$, and its maximal ideal
is generated by the regular element $\pi_i$.  Tensoring the exact
sequence
\[
 0\longrightarrow R\longrightarrow V\oplus k
 \xrightarrow{(v,c)\mapsto\overline v-c}\prod_i k_i
 \longrightarrow0
\]
with the flat algebra $R'$ identifies $R'$ with the fiber product of
these DVR factors over their common residue field $l$.  Theorem~\ref{branch-split}
gives $\K$ for $R'$, and faithful flatness descends it to $R$.

For necessity in residue characteristic two, suppose that $R$ is
Cohen--Macaulay, satisfies $\K$, and that two is either zero or a
nonzerodivisor.  If $R$ were not weakly normal, Lemma~\ref{wn}(3) would
give $q\in Q(R)\setminus R$ with $2q,q^2\in R$.  Theorem~\ref{realize} would
produce a nonsplit module-finite extension of finite projective dimension,
a contradiction.  Thus $R$ is weakly normal and in particular
seminormal.  Together with the sufficient direction this proves the
asserted equivalence under \eqref{separable-branches}.  The statement
about $F$-purity follows from the next corollary.
\end{proof}

The corresponding Frobenius consequence holds in every positive equal
characteristic.  It fits the classical structure theory of
one-dimensional $F$-pure rings \cite{GW,Gib}.

\begin{corollary}\label{separable-fpure}
Let $R$ be an excellent one-dimensional Cohen--Macaulay local ring of
characteristic $p>0$ satisfying \eqref{separable-branches}.  Then
seminormality, weak normality and $F$-purity are equivalent.
No finiteness hypothesis on Frobenius is required.
\end{corollary}
\begin{proof}
Weak normality implies seminormality.  If Frobenius is pure and
$q=b/a\in Q(R)$ satisfies $q^p\in R$, with $a$ a nonzerodivisor, then
$b^p\in a^pR$.  Tensoring Frobenius with $R/aR$ gives an injection
\[
 R/aR\longrightarrow F_*(R/a^pR),\qquad
 b+aR\longmapsto F_*(b^p+a^pR).
\]
Hence $b\in aR$, so $q\in R$, and Lemma~\ref{wn} gives weak normality.

Conversely, suppose that $R$ is seminormal.  The completion and finite
\'etale extension used above give a faithfully flat local extension
$R\to R'$ such that $R'$ is complete and is a fiber product of DVRs
with common residue field $l$.  A coefficient field of $R'$ maps to
a coefficient field of each branch.  Choosing uniformizers therefore
identifies
\[
 R'\simeq\prod_l l[[t_i]].
\]
The field $l$ need not be perfect.  Regard $F_*l$ as an $l$-vector
space with scalar action $a\cdot F_*c=F_*(a^pc)$.  Since $F_*1\ne0$,
choose an $l$-linear functional
\[
 \lambda:F_*l\longrightarrow l,\qquad \lambda(F_*1)=1.
\]
On each branch define
\[
 C_i\!\left(F_*\sum_{n\geq0}c_nt_i^n\right)
   =\sum_{m\geq0}\lambda(F_*c_{pm})t_i^m.
\]
Additivity and the identity
$\lambda(F_*(a^pc))=a\lambda(F_*c)$ show coefficientwise that $C_i$
is $l[[t_i]]$-linear for the Frobenius module structure.  It takes
$F_*1$ to $1$ and hence splits Frobenius.  The maps $C_i$ induce the
same functional $\lambda$ on the common residue field, so their
product restricts to a Frobenius splitting $F_*R'\to R'$.

Finally, $R\to R'\to F_*R'$ is pure: the first map is faithfully flat
and the second is split.  This composite factors through the
Frobenius map $R\to F_*R$.  After tensoring with any $R$-module,
injectivity of the composite forces injectivity of the first factor.
Thus $R$ is $F$-pure.
\end{proof}

For instance $k[[x_1,\ldots,x_s]]/(x_ix_j:i<j)$ identifies with $\prod_k k[[x_i]]$ by restriction to the axes. Theorem~\ref{branch-split} gives $\K$ for all $s$. For $s\ge2$ the normalization itself is a nonsplit module-finite extension because each local factor misses the other branches. Its projective dimension is infinite, since a single branch has generic ranks $1$ and $0$ at different minimal primes and hence cannot admit a finite free resolution.

\subsection{The complete Gorenstein curves}\label{forms}
\begin{corollary}\label{list}
Let $R$ be a complete one-dimensional Gorenstein local ring of residue characteristic $p>0$. Then $R$ is weakly normal precisely when one of the following holds.
\begin{enumerate}
\item $R$ is a complete DVR.
\item $R=V_1\times_k V_2$ with complete DVRs of residue field $k$.
\item $R$ is the inverse image of $k$ under $V\to k'$ for a complete DVR $V$ and a separable quadratic field extension $k'/k$.
\end{enumerate}
Under the hypotheses of the converse in Theorem A, these are exactly
the complete one-dimensional Gorenstein local rings satisfying $\K$.
\end{corollary}
\begin{proof}
Weak normality implies seminormality.  Apply
Lemma~\ref{conductor-description}, which requires no separability
assumption, and retain its notation $V,J,k_i$.  Set
\[
 d=\sum_i[k_i:k],\qquad \xi=(\pi_i)_i.
\]
The element $\xi$ is regular, $\xi J=J^2$, and
$\dim_kJ/J^2=d$.  Moreover, $\xi R/J^2$ is the line spanned by the
class of $\xi$.  Hence the maximal ideal of $R/\xi R$ has square
zero and dimension $d-1$.  If $d>1$, this maximal ideal is its socle,
and the Gorenstein property forces $d-1=1$.

For $d=1$, the inverse-image description gives the whole DVR.  For
$d=2$, there are either two branches with residue field $k$, or one
branch with a quadratic residue extension $k'/k$.  If this quadratic
extension were inseparable, its characteristic would be two and
there would be $\alpha\in k'\setminus k$ with $\alpha^2\in k$.
A lift $q\in V$ of $\alpha$ would satisfy $q\notin R$ and $q^2\in R$,
contradicting Lemma~\ref{wn}.  This gives the listed possibilities.

Conversely, in each displayed inverse image the ambient product of DVRs is finite over $R$ because lifts of a residue basis containing $1$ generate it over $R$, as $J$ already belongs to $R$. Eakin's theorem gives Noetherianity, and the componentwise powers of $J$ give completeness. The regular parameter $\xi$ and the socle computation above give the Gorenstein property in the two singular cases. Since $\xi V\subset R$, the total quotient ring is the product of the branch fraction fields. If $q^p\in R$, valuation on each branch puts $q$ in $V$. On the two rational branches equality of residue $p$th powers implies equality of residues. On a separable quadratic branch, an element whose $p$th power lies in $k$ already belongs to $k$. These observations prove root closure, including the DVR case.
\end{proof}

\begin{example}[Two mixed-characteristic nodes]
For every prime $p$, evaluation at $0$ and $p$ identifies
$\mathbf Z_p[z]/(z(z-p))$ with $\mathbf Z_p\times_{\mathbf F_p}\mathbf Z_p$ by $a+bz\mapsto(a,a+pb)$. This ring satisfies $\K$ although its reduction modulo $p$ is $\mathbf F_p[z]/(z^2)$, on which Frobenius has a kernel. The ring $\mathbf Z_p[[t]]/(pt)$ also satisfies $\K$, since evaluation at $t=0$ and reduction modulo $p$ identify it with $\mathbf Z_p\times_{\mathbf F_p}\mathbf F_p[[t]]$. The second example has a vertical branch on which $p$ vanishes.
\end{example}
\begin{example}[The boundary of the quadratic construction]
For $N\ge2$, evaluating at $0,p^N$ identifies $R_N=\mathbf Z_p[z]/(z(z-p^N))$ with the pairs congruent modulo $p^N$. The fraction $q=z/p=(0,p^{N-1})$ is missing, while $pq,q^p\in R_N$ because $p(N-1)\ge N$. For $p=2$ Theorem~\ref{realize} gives a nonsplit module-finite extension. For odd $p$ the determinant bound still restricts the rank of a possible counterexample, but the quadratic construction does not settle $\K$ for $R_N$. Likewise, necessity in mixed characteristic two when two is a nonzero zerodivisor is outside the scope of Theorem A.
\end{example}

\section{Splitting along the crossing locus}\label{geometry}
Let $k$ be any field of characteristic $p>0$. A \emph{splinter} is a ring over which every module-finite extension splits.
\subsection{Trace ideals and coefficient extraction}
Trace ideals are uniformly compatible with Cartier operators \cite{DT2023,PS2023}. The following formula works over an arbitrary coefficient field.
\begin{lemma}\label{cartier-trace}
Let $q=p^e$.  Every additive map $C:R\to R$ satisfying
$C(a^qb)=aC(b)$ preserves $J(S/R)$.
\end{lemma}
\begin{proof}
If $c=\phi(1)$, the map
\[
 S\longrightarrow R,\qquad s\longmapsto C\bigl(\phi(s^q)\bigr)
\]
is $R$-linear and takes $1$ to $C(c)$.  Additivity uses commutativity of
$S$, and scalar compatibility follows from
$\phi((as)^q)=a^q\phi(s^q)$.
\end{proof}

For $P=k[[X_1,\ldots,X_d]]$, choose a multi-index $\gamma$ with
$0\le\gamma_i<q$.  Choose also a $k^q$-linear functional $k\to k^q$ and
write $\lambda$ for its composite with inverse Frobenius $k^q\to k$.  Then
\[
 L_\gamma\!\left(\sum_\nu b_\nu X^\nu\right)
   =\sum_{\delta\in\mathbf N^d}
       \lambda(b_{q\delta+\gamma})X^\delta
\]
is a Cartier operator because coefficient convolution gives
$L_\gamma(A^qB)=A L_\gamma(B)$.
For $R=P/(f)$ it induces the operator
\[
 \overline H\longmapsto
 \overline{L_\gamma(f^{q-1}H)}.
\]
The map is well-defined because replacing $H$ by $H+fA$ changes the
output by $fL_\gamma(A)$.  This is the coefficient calculation underlying
Fedder's criterion \cite{Fed}.  If the coefficient of $X^\gamma$ in the
input is a nonzero element $b\in k$, extend $b$ to a $k^q$-basis of $k$
and choose the functional to send $b$ to $1$.  The extracted series then
has constant term $1$.  No assumption on $[k:k^p]$ is involved.

\begin{lemma}\label{conductor-lower}
Let $R$ be reduced and $R\subseteq B\subseteq Q(R)$ a module-finite birational extension. If $B$ is a splinter, then for every module-finite extension $R\subseteq S$ one has
\[
 (R:B)=J(B/R)\subseteq J(S/R).
\]
A subring which is a module direct summand of a reduced splinter is itself a splinter.
\end{lemma}
\begin{proof}
A functional on $B$ becomes scalar multiplication on $Q(R)$, so its values at $1$ form $(R:B)$. In $S\otimes_R Q(R)$, take the subring generated by $B$ and the image of $S$. It is module-finite over $B$ and hence has a $B$-linear retraction to $B$. Restrict to $S$ and multiply by any element of $(R:B)$ to obtain the desired functional $S\to R$.

For the second assertion let $B\subseteq H$ split as $B$-modules with $H$ a reduced splinter. For a module-finite extension $B\subseteq T$, the map $\Spec T\to\Spec B$ is surjective by lying over. It stays surjective after base change to $H$, so $H\to H\otimes_B T$ has nil kernel and is injective. Compose $T\to H\otimes_B T$ with the retractions to $H$ and to $B$.
\end{proof}

\subsection{Coordinate crossings and suspensions}
\begin{theorem}\label{cross-split}
Let $s\ge2$ and let
\[
 R=k[[x_1,\ldots,x_s,z_1,\ldots,z_n]]/(x_ix_j\mid i<j),
 \qquad P_0=(x_1,\ldots,x_s)R .
\]
Every evaluation ideal on $R$ is either $R$ or $P_0$, and the normalization realizes $P_0$. Splitting is detected at $P_0$, and $R$ has $\K$.
If $0\ne g\in(z_1,\ldots,z_n)^2$, the hypersurface
$k[[x,y,z_1,\ldots,z_n]]/(xy-g)$ is a domain and a splinter.
\end{theorem}
\begin{proof}[Proof for coordinate crossings]
Write $C=k[[z_1,\ldots,z_n]]$.
Represent $R$ as the tuples in $\prod_i C[[x_i]]$ having a common value
at $x_i=0$.  If $a\in J(S/R)\setminus P_0$, this common value has a
nonzero coefficient at some $z^\alpha$.  Choose $q>\max_j\alpha_j$
(any $q>1$ when $n=0$).  On each branch extract the residue exponents
$(0,\alpha)$ in $(x_i,z)$, using the same $\lambda$ on the coefficients.
The value at $x_i=0$ of each output depends only on the common input in
$C$.  The branch operators therefore preserve $R$ and define a Cartier
operator carrying $a$ to a unit.  Lemma~\ref{cartier-trace} implies
\[
 J(S/R)\ne R\quad\Longrightarrow\quad J(S/R)\subseteq P_0.
\]
The normalization is $B=\prod_i C[[x_i]]$. Indeed $B$ is a product of normal domains, is module-finite over $R$ by adjoining its coordinate idempotents, and agrees with $R$ after inverting $(x_i)_i$. Its conductor is $P_0$, since multiplying a tuple by each coordinate idempotent shows that its common value must vanish for the product to stay in $R$. Each factor is a regular splinter \cite{Andre}, so $B$ is a splinter by applying retractions factor by factor. Lemma~\ref{conductor-lower} yields $P_0\subseteq J(S/R)$ and $J(B/R)=P_0$. This proves the dichotomy.
Since evaluation ideals commute with localization, the extension splits
over $R$ if and only if it splits over $R_{P_0}$.

It remains to apply the curve theorem when $\pd_R S<\infty$.
Write $K=\operatorname{Frac}(C)$.  For $j\ge1$, the quotient $R/P_0^j$
is finite free over $C$, with basis $1$ and the pure powers $x_i^t$ for
$1\le t<j$.  Tensoring with $K$ gives the local Artinian truncations of
the union of axes over $K$.  Taking inverse limits gives
\[
 \widehat{R_{P_0}}\simeq\prod_K K[[x_i]].
\]
Localization and completion preserve the module-finite extension and its finite
projective dimension.  Theorem~\ref{branch-split} gives a retraction over this completed
curve.  Faithfully flat descent gives one over $R_{P_0}$, and the
localization criterion gives one over $R$.
\end{proof}

The strong $F$-regularity of the suspension $xy-g(z)$ is familiar in the
usual geometric setting \cite[p.~5]{HP2017}.  We retain the direct argument below to
cover arbitrary coefficient fields and all finite targets allowed in the
definition of $\K$.

\begin{proof}[Proof for suspensions]
Let $f=xy-g(z)$ with $0\ne g\in(z)^2$.  Replacing $xy$ by $g$ preserves
or increases total degree and gives a convergent normal form
\[
 c=\sum_{a\ge0}x^a A_a(z)+\sum_{b\ge1}y^b B_b(z).
\]
For uniqueness, substitute $y=g/x$ in $k((x))[[z]]$.  At each fixed
$z$-degree only finitely many negative powers of $x$ contribute, because
$g$ has positive order.  The coefficient of $x^{-b}$ is $g^bB_b$, while
that of $x^a$ is $A_a$.  Since $g\ne0$, the normal form is unique and the
quotient embeds in the domain $k((x))[[z]]$.

Let $c\ne0$ and suppose $A_a\ne0$.  Choose $\beta$ such that
$[z^\beta]g^aA_a\ne0$, and take $q$ larger than $a$ and all $\beta_i$.
The coefficient of $x^{q-1}y^{q-1-a}z^\beta$ in $f^{q-1}c$ is
\[
 (-1)^a\binom{q-1}{a}[z^\beta]g^aA_a
   =[z^\beta]g^aA_a\ne0.
\]
Indeed, the difference of the $x$- and $y$-exponents selects
$x^aA_a$, and
$(1+T)^{q-1}=\sum_{i=0}^{q-1}(-T)^i$ in characteristic $p$.
If only the second sum is nonzero, interchange $x$ and $y$.
Coefficient extraction therefore sends every nonzero element of $R$
to a unit by a suitable Cartier operator.

For any module-finite extension over a domain, $J(S/R)$ contains a nonzero
element: choose a linear functional on $S\otimes_R\operatorname{Frac}R$
taking $1$ to $1$ and clear its denominators on a finite generating set
of $S$.  Cartier compatibility now gives $J(S/R)=R$.
\end{proof}
\begin{corollary}\label{quadratic-rank}
Let $k$ have odd characteristic and write $X=(X_1,\ldots,X_d)$. If $f\in(X)^2$ has quadratic part of rank at least two, the ring $k[[X]]/(f)$ has $\K$.
It is a splinter if the quadratic rank is at least three, or if the residual series after splitting off a hyperbolic plane over a finite separable coefficient extension is nonzero.
\end{corollary}
\begin{proof}
In odd characteristic, diagonalize the quadratic form.  After a separable
extension of degree at most two, a nondegenerate binary summand becomes
hyperbolic.  Write the quadratic part as $uv+Q(z)$.
At each degree $m\ge3$, write the terms involving $u$ or $v$ as
$uA+vB$, with $A,B$ homogeneous of degree $m-1$.  The coordinate changes
$v\mapsto v-A$ and $u\mapsto u-B$ remove these terms without changing
the lower-degree terms.  Their successive compositions converge in the
maximal-ideal topology to a formal automorphism.  The equation becomes
$uv-g(z)$, with quadratic part of $g$ equal to $-Q$.
The crossing case applies if $g=0$ and the suspension case applies if $g\ne0$.
Quadratic rank at least three forces $Q\ne0$.  Finally, the finite
extension of coefficient fields commutes with power series and is
faithfully flat, so splitting descends by flat base change for the evaluation ideal.
\end{proof}

\subsection{Pinch points and the conductor}
\begin{theorem}\label{pinch}
Let $C=k[[z_1,\ldots,z_s,w_1,\ldots,w_n]]$, $s\ge1$, and $u\in C^\times$. Put $g=u z_1\cdots z_s$ and $R=C[[x,y]]/(x^2-y^2g)$.
If $p$ is odd, every evaluation ideal on $R$ is either $R$ or $(x,y)$. Its normalization realizes $(x,y)$, splitting is detected at $(x,y)$, and $R$ has $\K$. In characteristic two each of these hypersurfaces admits a nonsplit module-finite extension of generic rank four and projective dimension one.
\end{theorem}
\begin{proof}
First work in odd characteristic, and put
$C=k[[z_1,\ldots,z_n]]$.  More generally, let $0\ne g\in C$ satisfy
\begin{equation}\label{branchcond}
 \text{for each }0\ne c\in C\text{ some }q=p^e\text{ satisfies }
 c g^{(q-1)/2}\notin(z_1^q,\ldots,z_n^q).
\end{equation}
Write $R=C[[x,y]]/(x^2-y^2g)$ and $P_0=(x,y)R$.
Division by the monic equation in $x$ represents any element outside
$P_0$ as
$H=c+yA(y,z)+xB(y,z)$ with $0\ne c\in C$.
Choose $q$ as in \eqref{branchcond}, set $h=(q-1)/2$, and choose a
nonzero coefficient of $cg^h$ at $z^\beta$ with all $\beta_i<q$.
In $(x^2-y^2g)^{q-1}H$, the coefficient at
$x^{q-1}y^{q-1}z^\beta$ is precisely this coefficient.
The $xB$ term has odd $x$-degree, the $yA$ term overshoots the required
$y$-degree once the $x$-degree is fixed, and the contributing binomial
coefficient is
$(-1)^h\binom{q-1}{h}=1$.
Extraction sends $H$ to a unit.  Hence every proper $J(S/R)$ is
contained in $P_0$, and splitting is detected at $P_0$.

With $K=\operatorname{Frac}C$, truncation in $x,y$ gives
\[
 \widehat{R_{P_0}}\simeq K[[x,y]]/(x^2-gy^2).
\]
After the separable extension adjoining $\sqrt g$, if needed, the
coordinates $x+\sqrt g\,y$ and $x-\sqrt g\,y$ identify this ring with
a node.  Theorem~\ref{branch-split} therefore applies to every module-finite extension of finite
projective dimension.  Descent and the preceding detection criterion
give $\K$ for $R$.

For $g=u z_1\cdots z_s$ with $u\in C^\times$, choose a nonzero term
of $c$ of least total degree.  Its coefficient in $cu^h$ is multiplied
by the nonzero constant term of $u^h$, since $c$ has no terms of smaller
degree.  Multiplication by $(z_1\cdots z_s)^h$ translates its first
$s$ exponents by $h$.  All translated exponents are less than $q$
for sufficiently large $q=p^e$.  This proves \eqref{branchcond}, including the case of additional smooth variables.

For the ideal calculation, the coordinate change $z_1\mapsto u z_1$ absorbs the unit. Set $A=C[[y]]$ and $B=A[t]/(t^2-z_1\cdots z_s)$. The substitutions $z_i=v_i^2$ and $t=v_1\cdots v_s$ identify $B$ with the subring of
$H=k[[v_1,\ldots,v_s,y,w_1,\ldots,w_n]]$ consisting of series whose $v_i$-exponents are either all even or all odd. The two parity classes give the summands in $B=A\oplus At$. Projection onto these monomials is $B$-linear, since their parity vectors form a subgroup of $(\mathbf Z/2)^s$. Hence $B$ is a splinter by Lemma~\ref{conductor-lower}. For odd $p$, $B$ is the fixed ring of the sign changes $v_i\mapsto\epsilon_i v_i$ with $\prod_i\epsilon_i=1$. Indeed these fix precisely the two stated parity classes. Thus $B$ is normal, since $H$ is normal and an element of $\operatorname{Frac}B$ is fixed by every sign change.

The substitution $x=yt$ identifies $R=A\oplus yAt$ inside $B$. It is birational since $t=x/y$, and its conductor is $yB=(x,y)R$. To see this, $a+ybt\in R$, $a,b\in A$, belongs to the conductor exactly when its product with $t$ lies in $R$, or $a\in yA$. Thus $B$ is the normalization for odd $p$. Lemma~\ref{conductor-lower} supplies the lower bound $(x,y)\subseteq J(S/R)$, and the Cartier calculation gives the upper bound for every proper evaluation ideal. In particular $J(B/R)=(x,y)$, so $\pd_R B=\infty$ since $R$ has $\K$.

In characteristic two the same parity description and $R=A\oplus yAt$ give $t=x/y\notin R$ and $t^2\in R$. Apply Theorem~\ref{realize}. For the nonreduced ring $k[[x,y]]/(y^2)$, apply
Lemma~\ref{model} directly with $a=x$ and $\eta=y\notin xR$.
The resulting rank-two algebra has a minimal differential of size
$5\times3$, and the relation $x^6(y+xu)=0$ rules out a retraction.
\end{proof}

\section{Annihilators and the closed-fiber representation}\label{homol}
Let $R\subseteq S$ be a module-finite extension of a local ring with $\pd_R S\le1$. Fix the minimal resolution and notation
\begin{equation}\label{homol-resolution}
 0\longrightarrow R^b\xrightarrow{d}R^a
   \xrightarrow{\pi}S\longrightarrow0,
 \qquad H=\Ext^1_R(S,R),\qquad F=S/\m S.
\end{equation}
At each minimal prime of $R$, the localization of $S$ is free by
Auslander--Buchsbaum, since the localized base has depth zero.  Its rank
is $r=a-b$, and injectivity of $R\to S$ gives $r>0$.  Minimality therefore
gives
\[
 \dim_k(H/\m H)=b<a=\dim_k F.
\]
If $b=0$, the module $S$ is free and $1_S\notin\m S$ by Nakayama's lemma, so a coordinate functional gives a retraction.
We may consequently assume $b>0$ whenever we argue by contradiction.
Put $V=H/\m H$ and write $D_k(-)=\Hom_k(-,k)$.  Multiplication maps on
$S$ act on $H$ by contravariance. Commutativity makes this an
$S$-module action.  The induced actions on $V$ and on Tor factor
through $F$.

\begin{lemma}[The Tor obstruction]\label{torob}
If $R\subseteq S$ does not split, then
\[
 \operatorname{Tor}_1^R(H,k)\simeq D_kF
\]
as $F$-modules, and $V$ is not an $F$-module direct summand of this
Tor group.  Moreover,
\[
 \operatorname{Ann}_S H\subseteq J S,
 \qquad J=\{\phi(1_S):\phi\in\Hom_R(S,R)\}.
\]
\end{lemma}
\begin{proof}
Nonsplitting means $J\subseteq\m$.  In fact every $R$-linear functional
$\phi:S\to R$ then takes all its values in $\m$: for $s\in S$, the
functional $t\mapsto\phi(st)$ takes the value $\phi(s)$ at $1_S$.
Dualizing \eqref{homol-resolution} gives
\[
 0\longrightarrow S^*\longrightarrow(R^a)^*
   \xrightarrow{d^*}(R^b)^*\longrightarrow H\longrightarrow0.
\]
The entries of $d^*$ lie in $\m$, and the preceding observation gives
$S^*\subseteq\m(R^a)^*$.  Thus the two displayed free modules are the
first two modules of a minimal free resolution of $H$.  It follows that
\[
 \operatorname{Tor}_1^R(H,k)=(R^a)^*\otimes_R k
    \simeq D_k(S/\m S).
\]
A lift of multiplication by $s\in S$ to \eqref{homol-resolution}
acts on this dual presentation by transposition.  The identification
is therefore $F$-linear.

Suppose now that $V$ were an $F$-module summand of $D_kF$.
Duality would make $D_kV$ a summand of the regular $F$-module $F$,
so $D_kV\simeq eF$ for an idempotent $e\in F$.  In particular $e$
acts as the identity on both $D_kV$ and $V$.  Choose a lift $s\in S$
of $e$.  Multiplication by $s$ on $H$ is surjective by Nakayama's
lemma, because it induces the identity on $V$.  A surjective
endomorphism of a Noetherian module is an automorphism.  It follows
that $s$ acts invertibly on $\operatorname{Tor}_1^R(H,k)$.
Under the preceding identification this action is the action of $e$
on $D_kF$.  That action is faithful, and an invertible idempotent
endomorphism is the identity.  Hence $e=1$ in $F$, which would give
$\dim_kV=\dim_kF$, contrary to $b<a$.

Finally, let $s\in\operatorname{Ann}_S H$, and lift multiplication
by $s$ to matrices $P$ on $R^a$ and $Q$ on $R^b$, with
$Pd=dQ$.  Since $Q^*$ induces zero on $\operatorname{coker}d^*$,
its image is contained in $\operatorname{im}d^*$.  Freeness of
$(R^b)^*$ permits a lift $W:(R^b)^*\to(R^a)^*$ such that
$Q^*=d^*W$.  Setting $L=W^*$ gives $Q=Ld$.  Consequently
$(P-dL)d=0$, and there is a map $g:S\to R^a$ with
$P-dL=g\pi$.  Applying $\pi$ shows that $\pi g$ is multiplication
by $s$.  If $g_1,\ldots,g_a$ are the coordinate functionals of $g$,
then
\[
 s=\pi g(1_S)=\sum_{i=1}^a g_i(1_S)\pi(e_i)\in JS,
\]
as required.
\end{proof}

\subsection{An annihilator with nonzero linear term}
\begin{theorem}\label{linear-ann}
If $\operatorname{Ann}_R H\not\subseteq\m^2$, the extension splits.
On a one-dimensional Gorenstein base it suffices that $\m H^0_{\m}(S)=0$.
\end{theorem}
\begin{proof}
If $H=0$, minimality of the dual presentation forces $b=0$, and the
extension splits. Otherwise every annihilator of $H$ is a nonunit.
The hypothesis therefore supplies
$x\in\operatorname{Ann}_R H\cap(\m\setminus\m^2)$.
Choose a minimal free cover $P\to H$, with kernel $K$.  Then
\[
 xP\subseteq K\subseteq\m P,
 \qquad V=P/\m P,
 \qquad \operatorname{Tor}_1^R(H,k)=K/\m K.
\]
Multiplication by $x$ induces a well-defined map
\[
 \iota_x:V\longrightarrow K/\m K,
 \qquad \bar v\longmapsto xv\bmod\m K,
\]
because $x\m P\subseteq\m K$.  Choose a $k$-linear functional
$\lambda:\m/\m^2\to k$ with $\lambda(\bar x)=1$.
The inclusion $K\subseteq\m P$ gives a map
\[
 K/\m K\longrightarrow\m P/\m^2P
       \simeq(\m/\m^2)\otimes_k V.
\]
Composing it with $\lambda\otimes1_V$ gives a left inverse of
$\iota_x$.

Both maps are $F$-linear.  Indeed a lift to $P$ of multiplication
by an element of $S$ preserves $K$ and commutes with the scalar $x$. On $\m P/\m^2P$ it acts as the identity on the first tensor factor
and by its action on $V$ on the second.  Thus $V$ would be a direct
summand of $\operatorname{Tor}_1^R(H,k)$ in the nonsplit case,
contradicting Lemma~\ref{torob}.

For the stated consequence on a one-dimensional Gorenstein base,
set
\[
 T=H^0_{\m}(S)
   =\{s\in S:\m^n s=0\text{ for some }n\geq0\}.
\]
This is a finite-length module.  After completion, local duality
identifies $\Ext^1_{\widehat R}(S\otimes_R\widehat R,\widehat R)$
with the Matlis dual of $T$.  Here the Matlis dual is
$\Hom_{\widehat R}(T,E_{\widehat R}(k))$, with $E_{\widehat R}(k)$
the injective hull of the residue field.  Matlis duality is faithful
on finite-length modules, so for every $h\geq1$
\begin{equation}\label{torsion-annihilator-equivalence}
 \m^hT=0\quad\Longleftrightarrow\quad\m^hH=0.
\end{equation}
Ext commutes with this flat base change, and faithful flatness gives
the assertion over $R$.  If $\m T=0$, any minimal generator of
$\m$ is an annihilator with nonzero linear term.  A retraction over
the completion descends because its evaluation ideal is
$J\widehat R$, so $J\widehat R=\widehat R$ implies $J=R$.
\end{proof}

\subsection{The size of the closed-fiber action}
The numerical criterion can be stated without a restriction on the
number of algebra generators of the closed fiber.  What matters is
the dimension of the representation through which that fiber acts.

\begin{proposition}\label{representation-bound}
In the notation of \eqref{homol-resolution}, assume $b>0$ and
$\m^hH=0$ for some $h\geq1$.  Put
\[
 c=\dim_k\operatorname{im}
       \bigl(F\longrightarrow\operatorname{End}_k(V)\bigr).
\]
If
\begin{equation}\label{representation-inequality}
 a-c>b\,\ell_R(\m/\m^h),
\end{equation}
then $R\subseteq S$ splits.
\end{proposition}
\begin{proof}
The kernel of the representation has dimension $z=a-c$.
Since $V\ne0$, this kernel does not contain $1_F$.  Choose minimal
module generators $1_S,s_1,\ldots,s_{a-1}$ of $S$ so that the
classes of $s_1,\ldots,s_z$ form a basis of that kernel.
Write the rows of the matrix $d$ as
$\alpha,w_1,\ldots,w_{a-1}$, in this order, and set
\[
 N=R^{1\times b}/\langle w_1,\ldots,w_{a-1}\rangle,
 \qquad \eta=[\alpha]\in N.
\]
A functional taking $1_S$ to $1$ exists exactly when $\alpha$ lies
in the span of the other rows, or equivalently when $\eta=0$.
The dual presentation also identifies $H=N/R\eta$.  Hence
\begin{equation}\label{numerical-row-annihilator}
 \m^hN\subseteq R\eta.
\end{equation}

The rows of $d$ have entries in $\m$.  Consider the images of
$w_1,\ldots,w_z$ in $(\m/\m^h)^b$.  Their span has length at most
$b\ell_R(\m/\m^h)$, which is strictly smaller than $z$ by
\eqref{representation-inequality}.  In particular these $z$
generators cannot be minimal.  There is a relation among them
with a unit coefficient.  Complete its coefficient row to an
invertible matrix over $R$ and change the first $z$ rows accordingly.
After making the inverse change in the corresponding generators
of $S$, we may assume
\[
 w_1\in\m^h R^{1\times b}.
\]
The new $s_1$ still acts as zero on $V$, since its residue class is
a linear combination of the same kernel basis.

Lift multiplication by this $s_1$ to a pair of matrices $P,Q$ with
$Pd=dQ$.  The column of $P$ indexed by $1_S$ can be chosen to be
the basis vector indexed by $s_1$.  On the dual presentation of
$H$, the matrix $Q^*$ induces multiplication by $s_1$. Since
$d^*$ is minimal, its reduction modulo $\m$ is exactly the action
on $V$.  Thus $Q\equiv0\pmod\m$.  The row of $Pd=dQ$ indexed by
$s_1$ now reads
\[
 \alpha+\sum_{j=1}^{a-1}P_{1j}w_j=w_1Q.
\]
The right side belongs to $\m^{h+1}R^{1\times b}$.  Passing to
$N$ and using \eqref{numerical-row-annihilator}, we obtain
\[
 \eta\in\m^{h+1}N\subseteq\m R\eta.
\]
Nakayama's lemma applied to the cyclic module $R\eta$ gives
$\eta=0$, proving splitting.
\end{proof}

\begin{corollary}\label{schur-bound}
With the notation of Proposition~\ref{representation-bound}, put
$r=a-b$ and $d_h=\ell_R(\m/\m^h)$. Each of the following conditions suffices for splitting.
\begin{enumerate}[label=\textup{(\roman*)}]
\item $r>b d_h+\lfloor b^2/4\rfloor+1-b$.
\item $r>b d_h$ and $F$ is generated by at most two elements as a $k$-algebra.
\end{enumerate}
For $b\le3$, the bound $r>b d_h$ suffices for every closed fiber.
\end{corollary}
\begin{proof}
The image of $F$ is a commutative subalgebra of $\operatorname{End}_k(V)$, so the Schur--Jacobson bound gives
$c\le\lfloor b^2/4\rfloor+1$ \cite{Jac}. For two generators Gerstenhaber's theorem gives $c\le b$ \cite{Gerst,Berg}.
Both bounds hold over $k$, since scalar extension preserves the dimension of the image. Substitute in \eqref{representation-inequality}. For $b=1,2,3$ the first upper bound equals $b$.
\end{proof}

\subsection{Lifting a retraction from the tangent cone}\label{tangent}
\begin{theorem}\label{tangent-theorem}
Suppose $\operatorname{gr}_{\m}R$ is a geometrically integral splinter over $k$. If
$(\operatorname{Jac} S)^2 H=0$, the extension splits.
This applies in particular when
$R=\mathbf Z_p[[u,v,z,w_1,\ldots,w_d]]/(uv-z^2+h)$,
where $p$ is any prime and $h\in(p,u,v,z,w_1,\ldots,w_d)^3$.
\end{theorem}
\begin{proof}
First assume that $S$ is local, with maximal ideal $\mathfrak n$
and residue field $l$, and suppose the extension does not split.
We prove
\begin{equation}\label{strong-annihilator}
 \operatorname{Ann}_S H\subseteq\mathfrak n\m S.
\end{equation}
Lemma~\ref{torob} first gives
$\operatorname{Ann}_S H\subseteq JS\subseteq\m S$.
Choose minimal generators $x_1,\ldots,x_v$ of $\m$, and write an
annihilator as $u=\sum_i x_i a_i$ with $a_i\in S$.
Let $P\to H$ be a minimal free cover with kernel $K$, as above,
and lift multiplication by $a_i$ to $P_i:P\to P$.
Since $uH=0$, the operator $U=\sum_i x_iP_i$ takes $P$ into $K$
and induces a map
\[
 \bar U:V\longrightarrow K/\m K.
\]
This map is $F$-linear.  Indeed, if $P_s$ lifts multiplication by
$s\in S$, then $[P_s,P_i]$ maps $P$ into $K$, since the induced
actions commute on $H$.  Therefore
$[P_s,U]P\subseteq\m K$.
Contracting $K/\m K\to(\m/\m^2)\otimes_kV$ against the coordinate
functional dual to $x_i$ gives an $F$-linear map $\rho_i$ with
\[
 \rho_i\bar U=a_i|_V.
\]
If $a_i$ were a unit of $S$, its action on $V$ would be invertible. Composing $\rho_i$ with its inverse would split off $V$ from
$K/\m K$.  Lemma~\ref{torob} forbids this.  Thus each
$a_i\in\mathfrak n$, proving \eqref{strong-annihilator}.
The hypothesis $\mathfrak n^2H=0$ now gives
\begin{equation}\label{minimal-target-filtration}
 \mathfrak n^2=\m\mathfrak n,
 \qquad \mathfrak n^{j+1}=\m^j\mathfrak n\quad(j\geq1).
\end{equation}
Here the reverse inclusion $\m\mathfrak n\subseteq\mathfrak n^2$
is automatic, and the higher equalities follow by induction.

Put $P'=\operatorname{gr}_{\m}R$ and
$G=\operatorname{gr}_{\mathfrak n}S$.
Equation~\eqref{minimal-target-filtration} shows that $G$ is generated
in degrees zero and one as a finite $P'$-module.  It is also finite
over the graded algebra $l\otimes_kP'$, through the natural map to
$G$.  Geometric integrality makes $l\otimes_kP'$ a domain.  Since
$R\subseteq S$ is integral and $l/k$ is finite,
\[
 \dim(l\otimes_kP')=\dim R=\dim S=\dim G.
\]
A nonzero ideal in the affine domain $l\otimes_kP'$ has quotient
of smaller dimension, whereas finiteness gives
$\dim G=\dim((l\otimes_kP')/\ker)$.  Hence
\begin{equation}\label{graded-residue-inclusion}
 l\otimes_kP'\hookrightarrow G.
\end{equation}
In particular $P'\hookrightarrow G$ is a module-finite extension.
The splinter hypothesis gives a $P'$-linear retraction $G\to P'$.
Taking its component of matching degree on each homogeneous piece
produces a graded retraction $\gamma:G\to P'$.

Write $\lambda:l\to k$ for the degree-zero part of $\gamma$;
then $\lambda(1)=1$.  Choose a $k$-basis of $l$ beginning with $1$
and with the remaining basis elements in $\ker\lambda$.
Lift it to $\sigma_0=1_S,\sigma_1,\ldots,\sigma_{t-1}\in S$.
Their classes in $S/\m S$ are linearly independent, and we may
complete them to a basis by elements $y_1,\ldots,y_q\in\mathfrak n$.
On the corresponding free cover define an $R$-linear functional by
\[
 f(\sigma_0)=1,\qquad f(\sigma_i)=0\ (i>0),
 \qquad f(y_j)=L_j\in\m,
\]
where $L_j$ lifts $\gamma([y_j]_1)\in\m/\m^2$.

Consider any relation
\[
 \sum_i r_i\sigma_i+\sum_j a_jy_j=0.
\]
All its coefficients are in $\m$, by minimality of the generators.
In degree one of $G$, the second sum vanishes, since
$\m\mathfrak n=\mathfrak n^2$.  By
\eqref{graded-residue-inclusion} and independence of the residue
classes of the $\sigma_i$, we obtain $r_i\in\m^2$ for every $i$.
Taking degree two of the relation and applying the $P'$-linear
map $\gamma$ gives
\[
 r_0+\sum_j a_jL_j\in\m^3.
\]
Equivalently, in the notation of \eqref{homol-resolution},
$d^*f\in\m^3(R^b)^*$.
But $\m^2H=0$, so
\[
 \m^3(R^b)^*\subseteq\m\operatorname{im}d^*
                      =d^*(\m(R^a)^*).
\]
Choose $v\in\m(R^a)^*$ with $d^*v=d^*f$.
Then $f-v$ descends to $S$ and takes a unit value at $1_S$.
Normalizing that value gives a retraction, a contradiction.

For an arbitrary finite target $S$, pass to the completion of $R$.
The tangent cone is unchanged.  Since it is a domain, the completed
base is a domain because the initial forms of two nonzero elements have
nonzero product, and their orders are finite by Krull intersection.
The completed algebra is a product of finite local algebras.
At least one factor map from the completed base is injective.
Otherwise the product of the nonzero kernels would be a nonzero
ideal of that domain in the kernel of the map to the full product.  The Jacobson radical commutes
with this completion because in either algebra it is the inverse image
of the radical of the same finite-dimensional closed fiber.
Flat base change for Ext and passage to a direct factor preserve
both the projective-dimension hypothesis and the radical-square
annihilation.  Apply the local case to an injective factor, compose
with projection from the product, and descend using the evaluation
ideal as above.

Finally, consider the mixed-characteristic family in the statement.
The associated graded ring of the ambient regular local ring is
the polynomial ring
$\mathbf F_p[P,U,V,Z,W_1,\ldots,W_d]$.
The defining equation has initial form $UV-Z^2$.  Orders of
nonzero elements add in the ambient associated graded domain,
so the initial ideal of this principal ideal is generated by
$UV-Z^2$.  Thus the tangent cone is a polynomial extension of
\[
 k[U,V,Z]/(UV-Z^2)\simeq k[s^2,t^2,st],
 \qquad k=\mathbf F_p.
\]
This description remains a domain after every field extension.
Projection of $k[s,t]$ onto the monomials of even total degree is
linear over $k[s^2,t^2,st]$ and gives a retraction even when
$p=2$.  The same projection extends after adjoining the remaining
polynomial variables.  Polynomial rings over $k$ are splinters by
the coefficient-extraction argument of Section~\ref{geometry}.
The direct summand is a splinter by Lemma~\ref{conductor-lower}.
Finally $p$ is a nonzerodivisor on the stated hypersurface.
Indeed the defining equation remains nonzero modulo $p$ in the
ambient power series domain, so reduction modulo $p$ shows
$(uv-z^2+h):p=(uv-z^2+h)$.
All the asserted hypotheses are therefore satisfied.
\end{proof}

\section{Quadratic structure and a cubic multiplicity bound}\label{ranks}
For a torsion-free module $M$ of rank $r$ over a normal
domain, its \emph{reflexive determinant} is
\[
 \det M=\left(\bigwedge^r M\right)^{**},
 \qquad M^*=\Hom_R(M,R).
\]
If $R$ is local and $M$ has finite projective dimension, then
$\det M\simeq R$.  Indeed, choose a finite free resolution of $M$.
On the open set where $M$ is free, its determinant is the alternating
tensor product of the determinants of the free modules in this
resolution, hence is trivial.  That open set contains every point
of height at most one, since a torsion-free module over a height-one
localization of a normal domain is free over a discrete valuation ring.
An isomorphism of reflexive modules on this open set extends uniquely
across its complement, which has codimension at least two.

\subsection{A determinant refinement of quadratic splitting}

\begin{theorem}[Quadratic structure]\label{quadratic-order}
Let $R$ be a normal local domain satisfying $(R_2)$ and $(S_3)$, and
$R\subseteq S$ a torsion-free module-finite extension of rank two with
$\det S\simeq R$. Then $S^{**}=R\oplus R\theta$ is a quadratic algebra and
$S=R\oplus I\theta$ for an ideal with $I^{**}=R$. In particular the extension splits. Finite projective dimension of $S$ suffices for the determinant hypothesis.
\end{theorem}
\begin{proof}
Put $K=\operatorname{Frac}R$ and regard $S$ as a submodule
of its two-dimensional generic algebra $S_K=S\otimes_RK$.
Inside $S_K$ its reflexive hull is
\begin{equation}\label{quadratic-reflexive-intersection}
 T=S^{**}=\bigcap_{\operatorname{ht}\mathfrak p=1}S_{\mathfrak p}.
\end{equation}
Each member of this intersection is an algebra, so $T$ is an algebra
as well.  It is module-finite over $R$.  Consequently its elements
are integral over $R$, and normality gives $T\cap K=R$.
To check torsion-freeness of $T/R$, suppose $at\in R$ for a nonzero
$a\in R$ and $t\in T$. Then $t\in K\cap T=R$.

We next show that $T/R$ is reflexive.  At a prime of height at most
two, the base is regular and the reflexive module $T_{\mathfrak p}$
is free. If
$1_T\in\mathfrak pT_{\mathfrak p}$, multiplication would imply
$T_{\mathfrak p}=\mathfrak pT_{\mathfrak p}$, contradicting
Nakayama's lemma. Thus the unit is primitive and $(T/R)_{\mathfrak p}$ is free at these
primes.  At a prime of height at least three, $(S_3)$ gives
$\operatorname{depth}R_{\mathfrak p}\geq3$, while reflexivity of
$T$ gives $\operatorname{depth}T_{\mathfrak p}\geq2$.
The depth lemma applied to
\[
 0\longrightarrow R_{\mathfrak p}\longrightarrow T_{\mathfrak p}
   \longrightarrow(T/R)_{\mathfrak p}\longrightarrow0
\]
gives depth at least two for the quotient.  Hence $T/R$ is
torsion-free and satisfies $(S_2)$, which is the reflexivity
criterion over a normal domain.  This is precisely the quotient
step of \cite[Lemma~2.2]{HM}.

The map $t\mapsto1_T\wedge t$ induces a homomorphism
\[
 T/R\longrightarrow\det T=\det S.
\]
It is an isomorphism at every height-one prime, since there the
unit extends to a basis of the free rank-two module $T_{\mathfrak p}$.
Both sides are reflexive, so it is an isomorphism everywhere.
By the determinant hypothesis $T/R\simeq R$.  Lift a generator
of this quotient to $\theta\in T$.  Then
\begin{equation}\label{quadratic-order-description}
 T=R\oplus R\theta,
 \qquad \theta^2=\alpha+\beta\theta\quad(\alpha,\beta\in R).
\end{equation}
No division by two, separability assumption, or field hypothesis
on $S_K$ has been used.

For an element $a+b\theta\in S$, the scalar $a$ is already in $S$.
Subtracting it shows that
\[
 S=R\oplus I\theta,
 \qquad I=\{b\in R:b\theta\in S\}.
\]
The subset $I$ is an ideal.  Equation~\eqref{quadratic-reflexive-intersection}
gives $S_{\mathfrak p}=T_{\mathfrak p}$ at height one, so
$I_{\mathfrak p}=R_{\mathfrak p}$ there and $I^{**}=R$ as fractional
ideals.  Projection onto the scalar coordinate is the required
$R$-linear retraction.

The description also identifies the defect and the conductor:
\begin{equation}\label{quadratic-conductor}
 T/S\simeq R/I,
 \qquad (S:T)=\{t\in T:tT\subseteq S\}=IT.
\end{equation}
Indeed, writing $t=a+b\theta$, the condition $t\in S$ is $b\in I$,
and the condition $t\theta\in S$ is $a+\beta b\in I$.
Together they say $a,b\in I$, which is exactly $t\in IT$.
The finite-projective-dimension consequence follows from the
determinant observation at the beginning of this section.
\end{proof}

\subsection{Cubic algebras and Buchsbaum--Rim multiplicity}
In rank three the quotient by the unit has rank two, but need not
be free.  The proof measures its failure to be free by the length
of an Ext module.  The multiplicity input is the classical
parameter-module formula of Buchsbaum--Rim and its treatment in
\cite{BR,HH}.

\begin{theorem}[Cubic multiplicity bounds]\label{cubic-bounds}
Let $R$ be a three-dimensional Cohen--Macaulay normal local domain and
$R\subseteq S$ a torsion-free module-finite extension of rank three
with $\pd_R S\le1$. Put $E=\Ext^1_R(S/R,R)$,
$e=e_{\m}(R)$ and $v=\dim_k\m/\m^2$.
If $\ell_R(E)<4e$, or if $R$ is singular Gorenstein and
$v\ell_R(\m^2E)<4(e-1)$, the extension splits and $\mu_R(S)\le4$.
\end{theorem}
\begin{proof}
Put $C=S/R$.  This module is torsion-free.  Indeed, if $as\in R$
for $0\ne a\in R$ and $s\in S$, then in $S\otimes_R\operatorname{Frac}R$
the element $s$ is a scalar.  It is integral over $R$, hence belongs
to $R$ by normality.  The sequence
$0\to R\to S\to C\to0$ gives $\pd_R C\leq1$.

Under either numerical hypothesis
$E=\Ext^1_R(C,R)$ has finite length.  For the second hypothesis,
this follows because $\m^2E$ has finite length by assumption and
$E/\m^2E$ has finite length automatically.  Localizing a free
resolution of $C$ away from $\m$, the vanishing of $E$ makes the
dual presentation surjective and hence split.  Thus $C$ is free
on the punctured spectrum.  At the maximal ideal its depth is at
least two by Auslander--Buchsbaum.  Therefore $C$ is reflexive,
of rank two, and its determinant is trivial by finite projective
dimension.  Exterior product, followed by a choice of this
trivialization, gives
\begin{equation}\label{cubic-quotient-self-dual}
 C\simeq C^*.
\end{equation}
To justify this last isomorphism, the exterior pairing is perfect
at height-one primes and both modules are reflexive. The resulting
codimension-one isomorphism extends over $R$.

Let $b=\dim_k(E/\m E)$.  If $b>0$, the minimal resolution of $C$
and its dual have the forms
\begin{equation}\label{cubic-quotient-resolutions}
 \begin{gathered}
 0\longrightarrow R^b\longrightarrow R^{b+2}
       \longrightarrow C\longrightarrow0,\\
 0\longrightarrow C^*\longrightarrow R^{b+2}
       \longrightarrow R^b\longrightarrow E\longrightarrow0.
 \end{gathered}
\end{equation}
Let $N$ be the image in $R^b$.  It has finite colength and is
contained in $\m R^b$.  If $q=\mu_R(N)$, then $q\leq b+2$.
On the other hand, the ideal of maximal minors of a minimal
$b\times q$ matrix for $N$ is $\m$-primary, because $R^b/N=E$
is nonzero and of finite length.  The determinantal height bound
therefore gives $3\leq q-b+1$, so $q=b+2$.
This is exactly the number required for $N$ to be a parameter
module in a free module of rank $b$ over a three-dimensional base.

For such a parameter module over a Cohen--Macaulay ring, the
Buchsbaum--Rim formula gives
$\ell_R(R^b/N)=e_{\mathrm{BR}}(N,R^b)$ \cite{BR,HH}.
Here $e_{\mathrm{BR}}$ is the multiplicity determined by the leading
coefficient of
\[
 j\longmapsto
 \ell_R\bigl(\operatorname{Sym}^j(R^b)/N^j\bigr),
\]
where $N^j$ is the degree-$j$ image of the symmetric algebra of
$N$ in $\operatorname{Sym}(R^b)$: its leading term is
$e_{\mathrm{BR}}j^{b+2}/(b+2)!$.
The inclusion $N\subseteq\m R^b$ compares these colengths in every
degree and yields
\begin{equation}\label{cubic-BR-bound}
 \ell_R(E)=e_{\mathrm{BR}}(N,R^b)
     \geq e_{\mathrm{BR}}(\m R^b,R^b)
     =e\binom{b+2}{3}.
\end{equation}
For the last equality, the degree-$j$ colength for $\m R^b$ is
\[
 \binom{j+b-1}{b-1}\ell_R(R/\m^j).
\]
Since $\ell_R(R/\m^j)$ has leading term $ej^3/3!$, comparison
of the leading coefficients proves the equality.  As
$\binom{4}{3}=4$, the inequality $\ell_R(E)<4e$ forces $b\leq1$.

For the second criterion assume that $R$ is Gorenstein.
Two dimension shifts in \eqref{cubic-quotient-resolutions}, followed
by \eqref{cubic-quotient-self-dual}, give
\[
 \Ext^3_R(E,R)\simeq\Ext^1_R(C^*,R)\simeq E.
\]
Local duality identifies the first module with the Matlis dual
of $E$.  Completion preserves all the finite
lengths and the displayed isomorphisms.
Thus $E$ is Matlis self-dual, and its socle
$(0:_E\m)$ has length $\dim_k(E/\m E)=b$.
If $x_1,\ldots,x_v$ minimally generate $\m$, the map
\[
 \m E\longrightarrow(\m^2E)^v,
 \qquad u\longmapsto(x_1u,\ldots,x_vu),
\]
has kernel contained in that socle.  It follows that
\begin{equation}\label{cubic-socle-bound}
 \ell_R(E)\leq2b+v\ell_R(\m^2E).
\end{equation}
For $b\geq2$ one has $\binom{b+2}{3}\geq2b$.
Combining \eqref{cubic-BR-bound} and \eqref{cubic-socle-bound}
therefore gives
\[
 v\ell_R(\m^2E)\geq2b(e-1)\geq4(e-1).
\]
The last inequality uses $e\geq2$, since the Cohen--Macaulay base
is singular.  This contradicts the second hypothesis, so again
$b\leq1$.

If $b=0$, the quotient $C$ is free and the exact sequence
$0\to R\to S\to C\to0$ splits.  Suppose $b=1$.
Choose lifts $s_1,s_2,s_3\in S$ of the three minimal generators
of $C$.  If $(d_1,d_2,d_3)^t$ generates its relation module,
then $\sum_{i=1}^3d_i s_i\in R$.  Write this scalar as $-d_0$.
The generators $1_S,s_1,s_2,s_3$ present $S$ with the single
relation
\[
 d=(d_0,d_1,d_2,d_3)^t.
\]
Indeed any relation projects to a multiple of the given relation
of $C$, and subtraction leaves a scalar relation, which is zero
by injectivity of $R\subseteq S$.

Lift multiplication by $s_1$ to a matrix $P$ on this free cover,
choosing its first column to be the basis vector corresponding
to $s_1$.  Preservation of the relation module gives $Pd=\lambda d$
for some $\lambda\in R$.  Its row indexed by $s_1$ says
\[
 d_0+P_{11}d_1+P_{12}d_2+P_{13}d_3=\lambda d_1.
\]
Consequently $d_0\in(d_1,d_2,d_3)$.  A row vector with first
entry $1$ can therefore be chosen to kill $d$ and descends to
an $R$-linear retraction of the unit.  This last step is where
the algebra structure, beyond the rank-two module $C$, is used.
Finally $\mu_R(S)\leq1+\mu_R(C)=b+3\leq4$ under either bound.
\end{proof}

\section{Cubic extensions on the blow-up}\label{blowup}

Let $I\subseteq\m$ be as in Theorem C, or let $I=\m$ when $\operatorname{gr}_{\m}R$ is a hypersurface in four variables. In the latter case $R$ is automatically Cohen--Macaulay of dimension three. Indeed a homogeneous system of parameters in the graded hypersurface is a regular sequence. Its lifts are regular in $R$ by induction, using $\operatorname{gr}(R/aR)=\operatorname{gr}(R)/(a^*)$ whenever the initial form $a^*$ is a nonzerodivisor. Write
\[
 X=\operatorname{Bl}_I(\Spec R),\quad D=V(I\mathcal O_X),\quad
 L=I\mathcal O_X=\mathcal O_X(-D),\quad T=X\setminus D.
\]
The exceptional divisor is effective Cartier, $L$ is relatively ample, and $L|_D=\mathcal O_D(1)$ \cite{Stacks}. We identify $T$ with $\Spec R\setminus V(I)$. Since $I$ contains a regular sequence of length three, the local cohomology sequence for $V(I)$ gives $\Gamma(T,\mathcal O_T)=R$. We write $L^n\mathcal E$ for $L^{\otimes n}\otimes\mathcal E$, allowing negative $n$.

\subsection{Multiplication and the obstruction to splitting}
Let $U$ be the matrix presenting $S/R$ and write $s=b+2$. Dividing each column by a local equation for $D$ gives a bundle sequence
\begin{equation}\label{bundlec}
 0\longrightarrow (L^{-1})^b\longrightarrow\mathcal O_X^s
       \xrightarrow{\rho} V\longrightarrow0 .
\end{equation}
Indeed, the divided matrix on the closed fiber is $\overline U_I$. The rank-drop locus is closed in $X$ and proper over the local base. Its image misses the closed point, so it is empty. Thus \eqref{bundlec} is locally split. The bundle $V$ has rank two, $\det V=L^b$, and $V|_T=(S/R)|_T$. The kernel of $U:R^b\to R^s$ is supported in $V(I)$ and lies in a free module. A regular element of $I$ therefore forces the kernel to vanish. In particular $S/R$ has projective dimension one and $\pd_R S\le1$.

The algebra $S|_T$ is locally free of rank three. Formula~\eqref{perfect-trace} settles the case $3\in R^\times$. Henceforth the residue characteristic is three, so $2$ is a unit. Define
\[
 M=\mathcal Hom(\Sym^2 V,V),\qquad Q=\Sym^3(V^\vee)\otimes\det V.
\]
There is an exact sequence
\begin{equation}\label{cubicsequence}
 0\longrightarrow V^\vee\xrightarrow{i}M\xrightarrow{j}Q\longrightarrow0,
 \qquad j(\mu)(v)=v\wedge\mu(v,v),
\end{equation}
where $i(\ell)(v,w)=\ell(v)w+\ell(w)v$. In a basis write
$\mu_{11}=(a,b'),\ \mu_{12}=(c,d),\ \mu_{22}=(h,f)$.
The four coefficients of $j(\mu)$ are $b',2d-a,f-2c,-h$. They give surjectivity and identify the kernel with $i(V^\vee)$.

A local linear section $\sigma$ of $S|_T\to V|_T$ gives
$\mu_\sigma(v,w)=\overline{\sigma(v)\sigma(w)}$. Changing $\sigma$ by a linear form $\ell$ changes $\mu_\sigma$ by $i(\ell)$. Hence $j(\mu_\sigma)$ is a global binary cubic on $T$ \cite[Theorem~2.1]{Wood}.

\begin{lemma}\label{cubic-lifting}
This cubic extends to $\Phi\in H^0(X,Q)$. If $\Phi$ lifts to $H^0(X,M)$, then $R\subseteq S$ splits.
\end{lemma}
\begin{proof}
Lift the $s$ generators of $S/R$ to $S$ to obtain a map $h:\mathcal O_X^s\to S\otimes_R\mathcal O_X$. The quotient $S\otimes_R\mathcal O_X\to V$ is well defined since the image of $U$ is contained in that of the divided matrix. Denote it by $g$. The polynomial
\[
 z\longmapsto\rho(z)\wedge g(h(z)^2)
\]
on $\mathcal O_X^s$ restricts on $T$ to the pullback of the binary cubic. Locally split \eqref{bundlec}. All coefficients containing a kernel variable vanish on $T$, as polynomial coefficients over the local coordinate ring. A local equation of $D$ is a nonzerodivisor, so these coefficients vanish on $X$ as well. The remaining coefficients give the desired section, uniquely.

If $j(\mu)=\Phi$, write $\mu_{\sigma_i}-\mu=i(\ell_i)$ on opens of $T$ with local sections $\sigma_i$. On overlaps, injectivity of $i$ gives $\sigma_i-\sigma_j=\ell_i-\ell_j$. The sections $\sigma_i-\ell_i$ glue and split the unit sequence on $T$. Evaluate the resulting projection $S|_T\to\mathcal O_T$ on elements of $S$ and use $\Gamma(T,\mathcal O_T)=R$.
\end{proof}

Rank-two duality gives $L^{2b}Q=\Sym^3 V$ and
$L^{2b}M\simeq V\otimes\Sym^2 V$. The latter identification uses $2\in R^\times$ to polarize quadratic forms. Under these identifications $L^{2b}j$ is symmetric multiplication up to a sign. Thus a cubic monomial in the global generators of $V$ has a global lift through this twisted map.

\subsection{A vanishing in residue characteristic three}
\begin{lemma}\label{frobcube}
Let $Y$ be a scheme over a field of characteristic three with an invertible sheaf $\mathcal O_Y(1)$. Suppose
$H^0(Y,\mathcal O_Y(t))=0$ for $t<0$ and
$H^1(Y,\mathcal O_Y(t))=0$ for every $t$. If
\[
 0\longrightarrow\mathcal O_Y(-1)^b\longrightarrow
 \mathcal O_Y^{b+2}\longrightarrow E\longrightarrow0
\]
is a bundle sequence, then $H^0(Y,(\Sym^3 E)(t))=0$ for $t<-b$.
\end{lemma}
\begin{proof}
For a rank-two bundle in characteristic three the absolute Frobenius $F_Y$ gives an exact sequence
\begin{equation}\label{frob-cube-sequence}
 0\longrightarrow F_Y^*E\longrightarrow\Sym^3 E
 \longrightarrow E\otimes\det E\longrightarrow0 .
\end{equation}
In a basis $x,y$ the first arrow has image spanned by $x^3,y^3$. To see that the quotient is canonical, take the degree-three derivative
$\Sym^3 E\to E\otimes\Sym^2 E$, defined by the Leibniz rule. Its image is the image of the injection
\[
 E\otimes\det E\longrightarrow E\otimes\Sym^2 E,\qquad
 a\otimes(u\wedge v)\longmapsto u\otimes(av)-v\otimes(au).
\]
Indeed, the derivatives of $x^2y,xy^2$ give its two basis vectors up to sign, and the derivatives of the cubes vanish. This coefficient calculation works on nonreduced schemes.

The given resolution gives $H^0(Y,E(j))=0$ for $j<0$. Its Frobenius pullback is exact by local splitting and replaces the left-hand twist by $-3$. Thus $H^0(Y,(F_Y^*E)(t))=0$ for $t<0$ as well. Since $\det E=\mathcal O_Y(b)$, the two outer terms of \eqref{frob-cube-sequence} twisted by $t$ have zero global sections for $t<-b$.
\end{proof}

\begin{proof}[Proof of Theorem C]
In case (i) set $B=R/I$. Since $I$ is generated by a regular sequence, the conormal module is free of rank three and
\[
 \operatorname{gr}_I R=\Sym_B(I/I^2)\simeq B[X_0,X_1,X_2].
\]
Thus $D=\PP^2_B$ and its closed fiber is $Y=\PP^2_k$. In case (ii), $D=Y=V(f_2)\subset\PP^3_k$. Both closed fibers satisfy the line-bundle hypotheses of Lemma~\ref{frobcube}. For the quadric this follows by twisting $0\to\mathcal O_{\PP^3}(-2)\xrightarrow{f_2}\mathcal O_{\PP^3}\to\mathcal O_D\to0$. No smoothness hypothesis enters.

Let $E=V|_Y$. Restriction of \eqref{bundlec} gives precisely the presentation in Lemma~\ref{frobcube}, and
\[
 (L^n Q)|_Y\simeq (\Sym^3 E)(n-2b).
\]
For $0\le n<b$, the lemma gives zero global sections on $Y$. In case (i) the same vanishing holds on all of $D$. Here is the argument when $B$ is arbitrary. Put $\mathfrak a=\m/I$ and $\mathcal H=(L^nQ)|_D$. This bundle is flat over $B$, so its $j$th graded piece for the $\mathfrak a$-adic filtration is
$(\mathfrak a^j/\mathfrak a^{j+1})\otimes_k\mathcal H|_Y$. None of these pieces has a nonzero global section. A global section of $\mathcal H$ therefore belongs to $\mathfrak a^j\mathcal H$ for every $j$. Krull intersection makes its germ zero at each point over the closed point of $B$. Its support is proper over $B$, hence empty. This also explains why characteristic three is needed only on the closed fiber.

Use the restriction sequences for the Cartier divisor $D$ successively. The vanishing just proved places $\Phi$ in $H^0(X,L^bQ)$. To lift it, note that in both cases
\[
 H^1(D,\mathcal O_D(j))=0\quad(j\ge0),\qquad
 H^2(D,\mathcal O_D(j))=0\quad(j\ge-1).
\]
For $\PP^2_B$ these are the cohomology formulas for projective space over a ring \cite{Stacks}. For the quadric they follow from the same hypersurface sequence as above. Since the base is affine, relative Serre vanishing followed by descending induction in
\[
 0\longrightarrow L^{j+1}\longrightarrow L^j
       \longrightarrow\mathcal O_D(j)\longrightarrow0
\]
gives $H^1(X,\mathcal O_X)=H^2(X,L^{-1})=0$. Sequence \eqref{bundlec} now gives $H^1(X,V)=0$. Finally $L^b V^\vee\simeq V$, so \eqref{cubicsequence} twisted by $L^b$ lifts $\Phi$ to $H^0(X,L^bM)$. Lemma~\ref{cubic-lifting} gives the retraction. The closed-fiber matrix depends only on $U$ modulo $\m I$.
\end{proof}

\begin{corollary}\label{toeplitz}
Let $I=(a_0,a_1,a_2)$ be generated by a regular sequence in a local ring. For each $b\ge1$, take the $(b+2)\times b$ matrix whose $j$th column has $a_0,a_1,a_2$ in successive positions $j,j+1,j+2$ and zeros elsewhere. Every module-finite extension with unit quotient presented by this matrix, or by any perturbation of it in $\m I$, splits.
\end{corollary}
\begin{proof}
At a geometric point $[x_0:x_1:x_2]$ of $\PP^2_k$, the matrix sends a polynomial of degree $<b$ to its product with $x_0+x_1 Z+x_2 Z^2$. This polynomial is nonzero, so multiplication is injective. Apply Theorem C(i).
\end{proof}

For instance take $R=\mathbf Z_3[[x,y,z]]/(3^e+x^e+y^e+z^e)$ for any $e\ge2$. The ideal $(3,x,y)$ is a parameter ideal since its quotient is $\mathbf F_3[[z]]/(z^e)$, and $R$ is Cohen--Macaulay. Thus Corollary~\ref{toeplitz} applies for every $b$, including $b\ge4$, at arbitrarily large multiplicity. More generally, for any fixed triple as above, full column rank on $\PP^2_k$ defines an open subset of the affine space of coefficient matrices because the rank-drop incidence is closed and its projection is proper. The displayed matrices give a $k$-point of this open over every residue field.

\section{Higher-degree hypersurface tangent cones}\label{higher-surfaces}
We now blow up $I=\m$ and suppose that
$G=\operatorname{gr}_{\m}R=P/(f_e)$, with $P=k[X_0,\ldots,X_3]$,
$f_e$ homogeneous of degree $e\ge3$, so $R$ is Cohen--Macaulay of dimension three by the preceding section.
With $X,D,L,V,M,Q$ as in Section~\ref{blowup}, put
$E=V|_D$, $W=\Sym^3 E$, and let $b$ denote the number of columns of $U$.
As before, the linear part of $U$ must have full column rank on $D$.
By the trace argument, it suffices to work in residue characteristic three.
On a degree-$e$ surface, $H^2(D,\mathcal O_D(-1))^\vee=G_{e-3}$.
This is the group that vanished in the proof of Theorem C. The following lifting bound and presentation calculations account for it.
\begin{lemma}\label{high-order-cubic-lifting}
Every section of $L^{e+b-2}Q$ has a global preimage in
$L^{e+b-2}M$.
\end{lemma}
\begin{proof}
The hypersurface sequence on $\mathbf P^3_k$ and adjunction give
\[
 H^1(D,\mathcal O_D(n))=0\quad(n\in\mathbf Z),\qquad
 \omega_D=\mathcal O_D(e-4),
\]
and $H^2(D,\mathcal O_D(n))=0$ for $n\ge e-3$.
Relative Serre vanishing for $L$, followed by descending induction in
$0\to L^{n+1}\to L^n\to\mathcal O_D(n)\to0$, gives
\[
 H^1(X,L^n)=0\ (n\ge0),\qquad H^2(X,L^n)=0\ (n\ge e-3).
\]
Twisting \eqref{bundlec} by $L^{e-2}$ now gives
$H^1(X,L^{e-2}V)=0$.
By $V^\vee\simeq V\otimes L^{-b}$, this is the first cohomology of
the kernel in \eqref{cubicsequence} twisted by $L^{e+b-2}$.
The required surjectivity follows from its cohomology sequence.
\end{proof}

\begin{corollary}\label{three-entries}
Suppose $R\subseteq S$ is a module-finite extension with $S/R\simeq\operatorname{coker}U$ and the ideal $I_1(U)$ of entries has three generators. If the linear matrix has full column rank on $D$, then $R\subseteq S$ splits for every $b$ and every $e\ge3$.
\end{corollary}
\begin{proof}
The initial entries span at least three dimensions in $P_1$. Otherwise their common zero locus would contain a line, which meets $D$ over an algebraic closure. Since the entry ideal has three generators, the span has dimension three. Choose generators $a_0,a_1,a_2$ whose linear parts form a basis, and complete these forms to coordinates $X_0,X_1,X_2,Z$ of $P_1$. The point $[0:0:0:1]$ lies outside $D$. Consequently the coefficient of $Z^e$ in $f_e$ is nonzero. It follows that the linear parts form a system of parameters in $G$, and $a_0,a_1,a_2$ form a regular sequence in the Cohen--Macaulay ring $R$. Indeed the associated graded ring of their quotient is a quotient of $G/(X_0,X_1,X_2)$ and has finite length.
Projection from that point gives a finite surjection $D\to\PP^2_k$, as $f_e$ is monic in $Z$ after scaling. The linear matrix on $D$ is the pullback of $\overline U_{I_1(U)}$ on this plane. Its full rank descends to the plane, and Theorem C(i) applies.
\end{proof}

\subsection{Two relations}
\begin{theorem}\label{hyp-two}
Let $b=2$ and suppose $S/R\simeq\operatorname{coker}U$. If $e=3$ or the entries of the $4\times2$ matrix $U$ generate $\m$, then the module-finite extension $R\subseteq S$ splits. The same holds for any perturbation in $\m^2$.
\end{theorem}
\begin{proof}
Write $u,v$ for the columns of the linear part of $U$.
Define the graded $P$-module
\[
 N=\operatorname{coker}\bigl(
 P(-1)^{20}\longrightarrow P^4,\quad
 (A,B)\longmapsto Au+Bv\bigr),
\]
where $A$ and $B$ range over symmetric $4\times4$ matrices, and put
$N_G=N\otimes_PG$.
We first prove
\begin{equation}\label{two-relation-vanishing}
 (N_G)_q=0\qquad(q\ge3)
\end{equation}
under the stated hypotheses.
The standard generators $\epsilon_j$ of $N$ satisfy
$u_i\epsilon_j=-u_j\epsilon_i$ and $u_i\epsilon_i=0$.
Consequently
\[
 u_iu_j\epsilon_l=-u_i u_l\epsilon_j=u_lu_j\epsilon_i=-u_ju_i\epsilon_l.
\]
Since $2$ is a unit, this gives
$I(u)^2N=I(v)^2N=0$.
If the entries of $U$ generate $\m$, their initial forms span $P_1$, so $P_+=I(u)+I(v)$ and $(P_+)^3N=0$.

It remains to treat $e=3$ when those initial forms do not span $P_1$.
Their span then has dimension exactly three. A span of dimension at
most two would give a projective common zero locus of positive
dimension meeting the hypersurface $D$, contradicting the rank
assumption.
Choose coordinates so that all entries lie in
$P_0=k[X_0,X_1,X_2]$ and put $Z=X_3$.
Their common zero is outside $D$, so the coefficient of $Z^e$ in
$f_e$ is nonzero.  Projection from that point gives a finite surjective
map $D\to\mathbf P^2_k$.
The two columns are therefore independent at every geometric point of
$\mathbf P^2$.

Let $N_0$ be the corresponding module over $P_0$.
We show that $(N_0)_1=0$.
Write the columns as linear maps $A_1,A_2:k^3\to k^4$.
By duality against symmetric matrices, a functional on $(N_0)_1$
is represented by a $3\times4$ matrix $B$ for which both $A_iB$
are alternating.
The map $A_1$ is injective, so choose bases with
\[
 A_1=\binom{I_3}{0},\qquad A_2=\binom{H}{w}.
\]
Alternation of $A_1B$ forces $B=(B_0\ 0)$, where $B_0$ is an
alternating $3\times3$ matrix.
Alternation of $A_2B$ then gives $HB_0$ alternating and $wB_0=0$.
If $B_0\ne0$, its rank is two.
Both $B_0$ and $HB_0$ kill $\ker B_0$. The space of alternating
forms on $k^3$ with this property is one-dimensional.
Hence $HB_0=\lambda B_0$.
For any nonzero $v\in\operatorname{im}B_0$, we obtain
$Hv=\lambda v$ and $wv=0$, so $A_2v=\lambda A_1v$, a contradiction.
Thus $B=0$ and $(N_0)_1=0$.

Since $N_0$ is generated in degree zero and has no relations there,
$N_0=k^4$ in that degree and is zero in higher degrees.
It follows that
$N_G\simeq(k[Z]/(Z^e))^4$, after rescaling the coefficient of $Z^e$.
For $e=3$, this proves \eqref{two-relation-vanishing}.

We next translate this graded calculation into a lifting statement.
The degree-three part of the Koszul resolution of the symmetric
quotient on $D$ is
\[
 0\to\mathcal O_D(-2)^4\to\mathcal O_D(-1)^{20}
 \to\mathcal O_D^{20}\to W\to0 .
\]
It is exact because the bundle presentation is locally split, so its
Koszul generators become two independent polynomial variables.
Set $B'=\ker(\mathcal O_D^{20}\to W)$.
Serre duality and the vanishing of $H^1(\mathcal O_D(n))$ give
\begin{equation}\label{two-relation-obstruction}
 H^1(D,B'(t))^\vee=(N_G)_{e-2-t}.
\end{equation}
Indeed, dualizing the map on $H^2$ of the two leftmost terms gives
$(A,B)\mapsto Au+Bv$, with coefficients in degree $e-3-t$.
For $t<0$, the sequence
$0\to B'(t)\to\mathcal O_D(t)^{20}\to W(t)\to0$ identifies
$H^0(D,W(t))$ with $H^1(D,B'(t))$.
For $t\ge0$, it identifies the latter group with the cokernel of
the map from degree-$t$ coefficients times the twenty cubic monomials.
Here $H^0(D,\mathcal O_D(t))=G_t$.

Since $L^nQ|_D=W(n-4)$, for $0\le n<4$ the dual of its section
space is $(N_G)_{e+2-n}$.
Successive restriction therefore puts $\Phi$ in $L^3Q$ when $e=3$,
and in $L^4Q$ when $e\ge4$.
For $e=3,4$, Lemma~\ref{high-order-cubic-lifting} already applies.
For $e\ge5$ and $0\le t\le e-5$, equations
\eqref{two-relation-vanishing} and \eqref{two-relation-obstruction}
show that $H^0(D,W(t))$ is generated by those monomials.
Given a remainder in $L^{4+t}Q$, express its restriction in that form,
lift the coefficients from $G_t=\m^t/\m^{t+1}$ to $\m^t$, and use the
global generators of $V$.
Symmetric multiplication gives a global preimage under $j$ of the
resulting section.  Subtraction raises the remainder to $L^{5+t}Q$.
After finitely many steps it lies in $L^eQ$ and lifts by
Lemma~\ref{high-order-cubic-lifting}.
Lemma~\ref{cubic-lifting} then gives the retraction.
\end{proof}

For an explicit example, take
$R=\mathbf Z_3[[x,y,z]]/(3^e+x^e+y^e+z^e)$, where
$e=4r$, $r\ge1$, and $3\nmid r$.
Use the columns $(3,x,y,z)^t$ and $(x,y,z,3)^t$.
Their linear forms span $P_1$, and their rank drops only at
$[1:\lambda:\lambda^2:\lambda^3]$ with $\lambda^4=1$.
At every such point, the initial Fermat equation has value $4=1$ in
the residue field.
Thus the criterion applies to every module-finite extension with this
presentation, and also to every presentation obtained by perturbing
its entries in order at least two.

\subsection{Three relations and the exceptional section}
We may use a hypersurface equation in a regular ambient ring in this proof. Completion preserves the evaluation ideal by flat base change, and preserves $G$ and the linear matrix. By the Cohen structure theorem, $\widehat R$ is a quotient of a regular local ring \cite{Stacks}. In the kernel choose elements whose linear parts form a basis of the kernel on cotangent spaces. They are part of a regular system of parameters, so quotienting by them gives a regular local ring $(A,\mathfrak n,k)$ of dimension four mapping onto $\widehat R$. Since $\widehat R$ has depth three, Auslander--Buchsbaum gives $\pd_A\widehat R=1$. The nonzero kernel of $A\to\widehat R$ is consequently free of rank one, say $(F)$. The initial ideal of $(F)$ is generated by $\operatorname{in}F$ in $\operatorname{gr}_{\mathfrak n} A=P$. Under the chosen tangent coordinates $\operatorname{in}F$ is a scalar multiple of $f_e$. We can therefore work with $R=A/(F)$ without imposing a presentation of the original local ring.

Write $u,v,w$ for the linear columns of the $5\times3$ matrix $U$.
For symmetric $5\times5$ matrices, define
\begin{align*}
 d(B_1,B_2,B_3)&=(B_2u-B_1v,B_3u-B_1w,B_3v-B_2w),\\
 \delta(\alpha,\beta,\gamma)&=w^t\alpha-v^t\beta+u^t\gamma .
\end{align*}
\begin{theorem}[Three relations]\label{three-relation-criterion}
Suppose the $5\times3$ linear matrix has rank three everywhere on $D$ and its entries span $P_1$. Write
$d_q:(P_{q-1})^{45}\to(P_q)^{15}$ for $d$ in degree $q$. If
\[
                   \operatorname{rank} d_4=469,
\]
every module-finite extension $R\subseteq S$ with $S/R\simeq\operatorname{coker}U$ splits. The smaller test $\operatorname{rank}d_2=130$ suffices as well. For each $D$ the admissible matrices form a geometrically nonempty dense open, with dense rational points when $k$ is infinite.
\end{theorem}
\begin{proof}
Consider the graded complex
\[
 P(-1)^{45}\xrightarrow{d}P^{15}\xrightarrow{\delta}P(1)
\]
with $P^{15}$ in cohomological degree zero. Symmetry gives $\delta d=0$. Write $N_P$ for its middle cohomology and $N_G$ for the middle cohomology after reduction modulo $f_e$. By a constant change of basis, $\delta$ is the row of four variables followed by eleven zeros. Its cokernel is $k(1)$ and its kernel is generated in degrees zero and one by the extra basis vectors and the Koszul relations. Thus $N_P$ is generated in those degrees, and
\[
 \dim_k(\ker\delta)_q=15\dim_k P_q-\dim_k P_{q+1}\quad(q\ge0).
\]
The dimensions are $130$ in degree two and $469$ in degree four. Surjectivity onto the kernel in degree $q\ge2$ forces $(N_P)_j=0$ for every $j\ge q$. In particular the degree-four test follows from the degree-two test. The exact sequence of complexes given by multiplication by $f_e$ yields
\begin{equation}\label{three-relation-obstruction}
 0\longrightarrow N_P/f_eN_P\longrightarrow N_G\longrightarrow k(1-e)\longrightarrow0 .
\end{equation}
Here $f_e$ kills $k(1)$. In degrees $q\ge4$ the only possible nonzero part of $N_G$ is one copy of $k$ when $q=e-1$.

\emph{The exceptional section.}
The degree-three symmetric quotient now has the exact resolution
\[
 0\to\mathcal O_D(-3)\to\mathcal O_D(-2)^{15}\to
 \mathcal O_D(-1)^{45}\to\mathcal O_D^{35}\to W\to0.
\]
Exactness again follows by locally completing the three independent
linear forms to polynomial coordinates.
Set $B'=\ker(\mathcal O_D^{35}\to W)$.
Break the resolution at the image of its second arrow.
The vanishing of $H^1(\mathcal O_D(t))$ and Serre duality on the two
resulting short exact sequences identify
\[
 H^1(D,B'(t))^\vee=(N_G)_{e-2-t}.
\]
the transposed maps on the $H^2$ groups are exactly $d$ and $\delta$.
As before, for $t<0$ this also computes $H^0(D,W(t))^\vee$.

Now $L^nQ|_D=W(n-6)$. For $0\le n<\min(5,e+1)$ its sections are dual to $(N_G)_{e+4-n}$. This index is at least four and greater than $e-1$, so the section space vanishes. Successive restriction puts $\Phi$ in $L^{\min(5,e+1)}Q$. For $e=3,4$, Lemma~\ref{high-order-cubic-lifting} and then Lemma~\ref{cubic-lifting} give the retraction. Assume $e\ge5$. The next restriction lies in $H^0(D,W(-1))$, a one-dimensional space by \eqref{three-relation-obstruction} in degree $e-1\ge4$.

\emph{Lifting the exceptional section.}
Let $Y$ be the blow-up of $\operatorname{Spec}A$ at its closed point,
with exceptional divisor $\Delta=\mathbf P^3_k$.
The surjection of Rees algebras embeds $X$ in $Y$, with
$\mathcal O_Y(-\Delta)|_X=L$ and scheme-theoretic intersection
$X\cap\Delta=D$.
Lift the entries of $U$ to $A$ and divide by the local equation of
$\Delta$.  This gives
$B_1'=\mathcal O_Y(\Delta)^3\to B_0'=\mathcal O_Y^5$.
Its restriction to $X$ is the subbundle in \eqref{bundlec}.

Let $K_i$ be the degree-$i$ part of the Koszul complex of the
corresponding linear forms in the symmetric algebra of $B_0'$.
Place $\operatorname{Sym}^iB_0'$ in degree zero, so that its term in
degree $-j$ is
$\bigwedge^jB_1'\otimes\operatorname{Sym}^{i-j}B_0'$.
Koszul multiplication gives
\[
 (K_1\otimes K_2)(\Delta)\longrightarrow K_3(\Delta).
\]
In degree $-3$, this is
\[
 B_1'\otimes\bigwedge\nolimits^2B_1'
 \longrightarrow\bigwedge\nolimits^3B_1',
\]
twisted by $\Delta$.
It is exterior multiplication from nine copies of
$\mathcal O_Y(4\Delta)$ to one.
A basis tensor maps to a generator with coefficient $1$, so the
map is surjective in characteristic three as well.

For $0\le j\le3$, the line bundle $\mathcal O_Y(j\Delta)$ has no
higher cohomology.
For $j=4$, its only positive cohomology is $H^3$, and restriction is
an isomorphism
\[
 H^3(Y,\mathcal O_Y(4\Delta))
   \simeq H^3(\Delta,\mathcal O_\Delta(-4))=k.
\]
These statements follow first for nonpositive multiples by relative
Serre vanishing and descending induction, and then for $j=1,\ldots,4$
from restriction to $\Delta$ and the cohomology of projective space.

Consider the hypercohomology spectral sequence for
$(K_1\otimes K_2)(\Delta)$ on $Y$.
Its only positive-cohomology term is $k^9$ in bidegree $(-3,3)$.
There is no incoming differential, and every possible outgoing
target either vanishes or lies outside the complex.
The edge map therefore surjects from degree-zero hypercohomology
onto $k^9$.
After restriction to $\Delta$, all terms are sums of
$\mathcal O_\Delta(-j)$ for $1\le j\le4$.
Consequently this $k^9$ is the entire degree-zero hypercohomology
on $\Delta$.
Naturality of the edge maps and the preceding restriction
isomorphism show that restriction from $Y$ onto this group is
surjective.
Multiplication maps it onto the corresponding one-dimensional
hypercohomology of $K_3|_\Delta(-1)$.

To identify the image on $D$, use the short exact sequence of
complexes
\[
 0\to K_3|_\Delta(-e-1)
   \xrightarrow{\ f_e\ }K_3|_\Delta(-1)
   \to K_3|_D(-1)\to0.
\]
The first complex has zero degree-zero hypercohomology.
Indeed, all its twists are negative, so only the top cohomology
of the leftmost term can contribute in total degree zero.
The kernel of the corresponding differential is dual to the
cokernel of
\[
 (P_{e-1})^{15}\xrightarrow{\ \delta\ }P_e,
\]
which is zero because the entries of $\delta$ span $P_1$.
Restriction therefore injects the one-dimensional group
$\mathbb H^0(\Delta,K_3|_\Delta(-1))$ into
$H^0(D,W(-1))$.
The dimension calculation above makes this map an isomorphism.

On $X$, the locally split sequence \eqref{bundlec} makes $K_i|_X$
a locally free resolution of $\operatorname{Sym}^iV$.
Thus the restrictions of the source and target of Koszul
multiplication resolve $L^5M$ and $L^5Q$, respectively.
Naturality along $D\subset X\subset Y$ and $D\subset\Delta\subset Y$
now gives the surjection
\begin{equation}\label{exceptional-cubic-lift}
 H^0(X,L^5M)\xrightarrow{\ j|_D\ }H^0(D,W(-1)).
\end{equation}

\emph{Completion of the lifting.}
Choose $\mu\in H^0(X,L^5M)$ whose image under
\eqref{exceptional-cubic-lift} is the restriction of $\Phi$.
The remainder $\Phi-j(\mu)$ lies in $L^6Q$.
For $0\le t\le e-6$, equation \eqref{three-relation-obstruction}
gives $(N_G)_{e-2-t}=0$, since $4\le e-2-t\le e-2$.
The monomial-lifting argument used for two relations therefore
raises the order to $e+1$.
When $e=5$, order six already suffices.
Lemma~\ref{high-order-cubic-lifting}, with $b=3$, lifts the final
remainder, and Lemma~\ref{cubic-lifting} gives the retraction.
\end{proof}

\begin{proof}[Geometric nonemptiness in Theorem~\ref{three-relation-criterion}]
In the affine space of $5\times3$ matrices of linear forms,
the condition that the entries span $P_1$ is open.
Full rank on $D$ is open by properness of the rank-drop incidence projection.
On the spanning locus, rank $469$ for $d_4$ is open by the minor criterion, and is the maximal possible rank since $\delta d=0$.
Let $\mathcal U_D$ be the intersection of these three opens.

To prove geometric nonemptiness, consider
\[
 U_0=\begin{pmatrix}a&0&0\\ b&a&0\\ c&b&a\\ d&c&b\\0&d&c\end{pmatrix}
\]
in four parameter coordinates.
The unique point where its rank drops is $[0:0:0:1]$. If the first nonzero coordinate
is among the first three, the corresponding three consecutive rows
give a nonzero triangular minor.
Its entries span $P_1$, and Appendix~\ref{cubcert} certifies
$\operatorname{rank}d_2=130$ over $\mathbf F_3$, hence over every
extension of that field.
Over $\bar k$, a linear coordinate change moves its unique bad point
outside the fixed hypersurface $D$.
Such a change preserves all graded ranks. The generation of $N_P$ in degrees zero and one then gives the degree-four test.
Thus $\mathcal U_D$ is geometrically nonempty and dense.

If $k$ is infinite, $\mathcal U_D(k)$ is dense by density of rational
points in affine space over an infinite field.
Every lift of a matrix in $\mathcal U_D(k)$ satisfies
the splitting criterion just proved.
All three conditions depend only on the linear part, so the same is
true after arbitrary perturbations of the entries in $\m^2$.
\end{proof}

For a family over the finite residue field, take
$A=\mathbf Z_3[[x,y,z]]$ with parameters $(3,x,y,z)$ and
$F=3^e+x^e+y^e+z^e$, for any $e\ge3$.
The initial equation takes value $1$ at the bad point of $U_0$, so
that matrix itself lies in $\mathcal U_D(\mathbf F_3)$.
If $3\nmid e$, the Fermat surface $D$ is smooth.
At a nonmaximal prime of $R$, at least one of $x,y,z$ is invertible:
otherwise $F=0$ would force $3$ into the prime as well.
The corresponding partial derivative of $F$ is a unit, proving that
$R$ is regular away from its closed point.
Since $R$ is Cohen--Macaulay, Serre's criterion makes it normal.
These isolated singularities have unbounded multiplicity, and
\[
 h^2(D,\mathcal O_D)=\binom{e-1}{3}.
\]
Thus the criterion does not rely on the vanishing of this exceptional
cohomology.  The exceptional section lifts by \eqref{exceptional-cubic-lift}.

For a general linear presentation on a higher-degree hypersurface, full column rank alone leaves the cubic lifting problem open in this argument.

\appendix
\section{Operator identities for the universal algebras}\label{cert}
Work with $P,x,Y$ as in Lemma~\ref{model}, and put $h=Y$ in equal characteristic and $h=2Y$ in mixed characteristic. In the ordered basis $1,u,v,w,d$ before passage to the quotient, set
\begingroup\small
\[
 U=\begin{pmatrix}
 0&x^2Y&0&x^2Y&0\\
 1&x^3+h&(1+x^2)Y&x^3+xY&x^2Y\\
 0&0&0&0&0\\
 0&0&x^2+(1+x)Y&0&(x+x^2)Y\\
 0&1+x&Y&0&0
 \end{pmatrix},
\]
\[
 V=\begin{pmatrix}0&0&0&0&0\\0&(1+x^2)Y&0&0&0\\
 1&0&0&0&0\\0&x^2+(1+x)Y&0&Y&xY\\0&Y&0&0&Y
 \end{pmatrix},
\]
\[
 W=\begin{pmatrix}0&x^2Y&0&0&0\\0&x^3+xY&0&0&(x^2+x^3)Y\\
 0&0&0&0&0\\1&0&Y&0&x^2Y\\0&0&0&0&xY
 \end{pmatrix}.
\]
\endgroup
We give a reproducible column test for preservation of the relations and for commutation on the quotient. For a column $z=(z_0,\ldots,z_4)^t$ form
\[
 B=z_1/x^4,\qquad A=(z_3-YB)/x^3,\qquad C=z_4/x^2
\]
initially in $P[1/x]$. Then $z\in\operatorname{im}D$ precisely when $A,B,C\in P$, $z_0=x^3YB$ and $z_2=YA$. Indeed these are the five equations $z=D(A,B,C)^t$, with three of them already solved. Applying this test to the columns of
\[
 UD,\quad VD,\quad WD,\quad UV-VU,\quad UW-WU,\quad VW-WV
\]
gives zero remainder in each case. For example, the first column of $VD$ has $z_1=z_4=0$, $z_3=x^3Y$, and $z_0=z_2=0$, so $A=Y$, $B=C=0$. The script \texttt{checks/verify.py} carries out the same test on all columns. Each division is exact in $\mathbf F_2[X,Y]/(Y^2)$ or $\mathbf Z[Y]/(Y^2)$, respectively.

\section{A rank certificate for the cubic presentation}\label{cubcert}
In the coordinates of $U_0$ in Section~\ref{higher-surfaces}, write its linear columns as $u,v,w$. Index the source of $d_2$ by a choice of one of three symmetric matrices, one of its fifteen independent entries, and a degree-one monomial. Index the target by a pair of columns, a row, and a degree-two monomial. In these bases $d_2$ is $150\times180$, while $\delta_2$ is $20\times150$. The matrices are built from the displayed formulas for $d,\delta$. An off-diagonal symmetric basis element contributes at both transposed positions. Working over $\mathbf F_3$, elimination gives
\[
 \operatorname{rank}\delta_2=20,\qquad \delta_2d_2=0,\qquad
 \det (d_2)_{I,J}=2
\]
for a $130\times130$ submatrix. The first two identities give the upper bound $130$ and the last attains it. The accompanying standard-library Python script \texttt{checks/verify\_cubic.py} constructs both matrices, prints the selected row and column indices and verifies the determinant on the original submatrix, separately from pivot selection.

\section*{Declaration of generative AI use}
Generative AI assisted with drafting, revising proofs, bibliographic comparisons, and preparing the symbolic verification scripts.

\end{document}